\documentclass{amsart}

\usepackage[top=3cm,bottom=3cm,inner=3.5cm,outer=3.5cm]{geometry}
\usepackage{tikz}
\usepackage{amssymb}
\usetikzlibrary{decorations.markings}
\usetikzlibrary{arrows}

\usepackage{comment}

\newcommand{\nc}{\newcommand}
\nc{\G}{{\Gamma}} \nc{\BC}{{\mathbb C}} \nc{\BQ}{{\mathbb Q}}
\nc{\BR}{{\mathbb R}} \nc{\BZ}{{\mathbb Z}} \nc{\BP}{{\mathbb P}}
\nc{\BN}{{\mathbb N}} \nc{\BM}{{\mathbb M}}
\nc{\fH}{{\mathbb H}}

\nc{\mat}{{\binom{a\,\ b}{c\,\ d}}}
\nc{\U}{{\mathcal U}}
\nc{\PS}{{\mbox{PSL}_2(\BZ)}} \nc{\SL}{{\mbox{SL}_2(\BZ)}}
\nc{\SR}{{\mbox{SL}_2(\BR)}} \nc{\PR}{{\mbox{PSL}_2(\BR)}}
\nc{\GL}{{\mbox{GL}}} \nc{\PQ}{{\mbox{PGL}_2^+(\BQ)}}
\nc{\GR}{{\mbox{GL}_2^+(\BR)}} \nc{\PG}{{\mbox{PGL}_2^+(\BR)}}
\nc{\GC}{{\mbox{GL}_2(\BC)}}
\nc{\f}{{\mathcal{F}(\fH)}}
\nc{\Cc}{\widehat{\BC}}
\nc{\e}{{E_{\rho}(\G)}}
\nc{\g}{{\gamma}}
\nc{\vm}{{V_{\rho}(\G)}}
\nc{\oo}{{\mathcal O}}
\nc{\M}{{\mbox{M}}}
\nc{\om}{{\omega}}
\nc{\Om}{{\Omega}}
\nc{\TX}{{\widetilde{X}}}
\nc{\ol}{\overline}
\nc{\cl}{{\mathcal L}}
\nc{\ce}{{\mathcal E}}
\nc{\la}{{\lambda}}
\nc{\La}{{\Lambda}}
\nc{\cz}{{\mathcal Z}}

\newtheorem{numbered}{}[section]
\newtheorem{thm}[numbered]{Theorem}

\newtheorem{remark}[numbered]{Remark}

\numberwithin{equation}{section}

\newcommand{\secref}[1]{\S\ref{#1}}

\newtheorem{exple}{\textbf{Example}}[section]

\begin{document}

\title[]{Modular groups and  planar  maps }
\author[]{Abdellah Sebbar} 
\author[]{Khalil Besrour}
\address{Department of Mathematics and Statistics, University of Ottawa, Ottawa Ontario K1N 6N5 Canada}
\email{asebbar@uottawa.ca}
\email{kbesr067@uottawa.ca}

\begin{abstract}
	In this paper we  give an explicit formula for the number of subgroups of the modular group of a given index that are genus zero and torsion-free and a formula for their conjugacy classes. We do so by exhibiting a correspondence between these groups and the  trivalent maps on a sphere.  We  focus on the particular case of  the subgroups of index 18 which have some interesting geometric properties.
\end{abstract}
\maketitle
\section{Introduction}
The purpose of this article is to enumerate the  torsion-free genus zero subgroups of  finite index of the modular group $\PS$ as well as their conjugacy classes. These enumerations will be given by means of a closed formula and a recursive one. In \cite{seb1}, the first author has  classified the same category of subgroups which are congruence. They form 33 conjugacy classes with indices ranging from 6 to 60. However, for non-congruence subgroups, there is no upper limit for the index, and it turns out that the number of such subgroups or their conjugacy classes increase rapidly with their index in $\PS$.

A subgroup $\G$ of finite index  $\mu$ of $\PS$ is characterized by its signature $(\mu;g,e_2,e_3,h)$ where $g$ is the genus of the modular curve
$X(\G)^*={\G}\backslash (\fH\cup \BQ\cup\{\infty\})$ where $\fH$ is the complex upper half-plane,  $e_r$ ($r=2,\,3$) is the number of inequivalent elliptic fixed points of order $r$, and $h$  is the number of inequivalent cusps of $\G$ satisfying the Riemann-Hurwitz formula
\[
g\,=\, 1+\frac{\mu}{12}-\frac{e_2}{4}-\frac{e_3}{3}-\frac{h}{2}.
\]
According to \cite{mil}, if $\mu$ and $h$ are positive integers and $e_2$,
 $e_3$ and $g$ are non-negative integers satisfying this equation, there exists a subgroup of $\PS$ with signature $(\mu;g,e_2,e_3,h)$. This correspondence works by associating to such data a class of pairs of permutations $[\phi,\psi]_{x_1}$ on a  set of $\mu$ letters such that $\phi$ fixes $e_2$ letters, $\psi$ fixes $e_3$ letters and $\phi\psi$ splits into $h$ disjoint cycles. The equivalence class is modulo a conjugation by a permutation $\lambda$ that fixes a chosen element $x_1$.  If $\G$ is a given subgroup of index $\mu$ of the modular group, then the set in question is the collection of cosets and the permutations $\phi$ and $\psi$ correspond respectively to the action on the cosets 
of the matrices 
$\displaystyle S=\binom{0\ \ -1}{1\ \  \ \ \ 0}$ and $\displaystyle V=\binom{1\  -1}{1\ \ \ \ 0}$. If one ignores the element $x_1$, then we obtain a correspondence between  the classes of pairs $[\phi,\psi]$ and the conjugacy classes of subgroups with signature  $(\mu;g,e_2,e_3,h)$.

When $e_2=e_3=g=0$, that is when the groups are genus zero and torsion-free, we have $\mu=6(h-2)$ and there is a way to attach to each group a unique trivalent planar graph on the sphere that we call a planar map. The vertices of such a map are the orbits of $\psi$ and so they are $\mu/3$ in number. Two vertices are joined by an edge if they are in the same $\phi-$orbit.

On the other hand, the theory of planar maps is  well developed by pioneering works by Tutte, Liskovets, Walsh, Wormald and others, and the census of these maps arose out  through enumeration techniques of the so-called rooted and unrooted maps.  It turns out that the genus zero and torsion-free subgroups of the modular group are in one-to-one correspondence with the rooted trivalent planar maps, while the conjugacy classes of these subgroups correspond to the unrooted ones. This allows us to derive enumeration formulas for such subgroups of a given index and for their conjugacy classes inside
$\PS$.

Furthermore,  as the torsion-free genus zero yield elliptic surfaces over ${\mathbb P}^1$, we make use of the classification of subgroups of index 18 to determine the $K3$ surfaces that possess an elliptic fibration  over ${\mathbb P}^1$ with 4 semi-stable singular fibers of type $I_n$  and a singular fiber of type $I_n^*$ following  Kodaira's notation, in the same way that the index 24 subgroups classify the semi-stable elliptic surfaces with six singular fibers.

This paper is organized as follows: In Section 2, we describe the structure of subgroups of the modular group and their connection with permutation groups. In Section 3, we explore the theory of rooted and unrooted maps and their enumeration and we provide several examples. In Section 4, we transform the Schreier graphs attached to a modular subgroup into rooted and unrooted maps on the sphere, yielding explicit formulas for the subgroups of a given index as well as their conjugacy classes. Finally, we describe briefly the case of the genus 1 subgroups that are torsion-free which correspond to smooth modular curves of genus 1. These are associated to rooted and unrooted maps on the torus.

\section{Modular Subgroups}\label{sec1}
Let $\SL$ be the modular group, that is the group of 2$\times$2 matrices with integral entries and determinant 1 and let $\PS=\SL/{\pm I}$. It acts on the upper half-plane $\fH$ by linear fractional transformations.
We refer to a finite index subgroup of $\PS$ as a modular subgroup. Such a subgroup $\G$ has the following presentation:

Generators: $A_1$, $B_1$, $\dots$, $A_g$,  $B_g$, $E_1$, $\dots$, $E_r$, $P_1$, $\cdots$, $P_h$.

Relations: $E_1^{m_1}=\dots =E_r^{m_r}=P_1\cdots P_h=A_1B_1A_1^{-1}B_1^{-1}\cdots A_gB_gA_g^{-1}B_g^{-1}=I$,

\noindent
where the generators $E_i$, $1\leq i\leq r$, are elliptic and $r$ is the number of conjugacy classes of cyclic
subgroups of ${\Gamma}$ of order $m_1,\dots m_r$, $m_i\in\{2,3\}$.
The generators $P_i$, $1\leq i\leq h$, are parabolic and $h$ is the number of inequivalent cusps, also referred to as the parabolic class number. Finally, the generators $A_i$, $ B_i$, $1\leq i\leq g$, are hyperbolic and $g$ is the genus of the compact Riemann surface $X(\G)={\G}\backslash (\fH\cup \BQ\cup\{\infty\})$.
Let $e_k$, $k\in\{2,3\}$, be the number of inequivalent elliptic fixed points of order $k$. The data  
\begin{equation}\label{sgn}
(\mu;g,e_2,e_3,h)
\end{equation}
 is referred to as the signature of $\Gamma$ and it is invariant under modular conjugacy of $\Gamma$. We also have the Riemann-Hurwitz formula \cite{ran}
\begin{equation}\label{r-h}
g\,=\, 1+\frac{\mu}{12}-\frac{e_2}{4}-\frac{e_3}{3}-\frac{h}{2}.
\end{equation}
According to \cite{mil}, for each integer $\mu>0$, $h>0$, $e_2\geq 0$, $e_3\geq 0$ and $g\geq 0$ satisfying \eqref{r-h}, there exists a modular subgroup $\G$ with signature $(\mu;g,e_2,e_3,h)$. Because it is relevant to this paper, we will make explicit how this correspondence works.

Let $X$ be a finite set of $\mu$ elements and $x_1\in X$ fixed. We consider pairs $(\phi,\psi)$ of permutations of $S_\mu$ acting on $X$ such that $\phi^2=\psi^3=I$ and $\Sigma=\langle\phi,\psi\rangle$ is transitive on $X$. Two such pairs $(\phi_1,\psi_1)$ and $(\phi_2,\psi_2)$ are said to be equivalent if there exists a permutation $\lambda\in S_\mu$ such 
\begin{equation}\label{equiv}
\lambda^{-1}\phi_1\lambda=\phi_2\,,\ \lambda^{-1}\psi_1\lambda=\psi_2\,, \ \lambda x_1=x_1.
\end{equation}
The equivalence class of a pair $(\phi,\psi)$ is denoted by $[\phi,\psi]_{x_1}$. There is a one-to-one correspondence  between modular subgroups $\G$ of finite index $\mu$ and classes  $[\phi,\psi]_{x_1}$, see ~\cite{mil}. The group $\G$ has signature $(\mu;g,e_2,e_3,h)$ and cusp decomposition
$\displaystyle \mu=\sum_{1\leq i\leq h}\,n_i$ if and only if
\begin{enumerate}
	\item $\phi$ fixes $e_2$ elements on $X$,
	\item $\psi$ fixes $e_3$ elements of $X$, and
	\item $\phi\psi$ splits into $h$ disjoint cycles of length $n_i$, $1\leq i\leq h$.
\end{enumerate}
Here $\Sigma$ is isomorphic to $\PS/G$ where $G$ is the largest normal subgroup of $\PS$ contained in $\G$ and $X$ is the set of  left cosets of $\G$ in $\SL$. The integers $n_i$ are the cusp widths which are defined as follows. Let $c\in\BQ\cup\{\infty\}$ be a cusp of $\G$ fixed by a parabolic element $P\in\G$, and let $\alpha\in\PS$ such that
$\alpha\, c=\infty$. If $T\,z=z+1$, then $\alpha P\alpha^{-1}=T^n$ as it fixes  $\infty$. The cusp width of $c$ is the smallest positive integer $n$ such that $\alpha^{-1}T^n\alpha\in\G$.

 If the fixed element $x_1$ is omitted in the above construction, then we obtain a one-to-one correspondence between equivalence classes  $[\phi,\psi]$ of pairs of permutations $(\phi,\psi)$ satisfying \eqref{equiv} (without having $\lambda$ to fix $x_1$) and modular subgroups up to modular conjugacy. 

Recall that a congruence subgroup $\G$ is a modular subgroup containing a principal congruence subgroup $\G(N)$ for some positive integer $N$, and the level of $\G$ is the smallest such $N$. It turns out that as the index grows, there are a lot more non-congruence groups than congruence ones. In fact the ratio of the two categories with the same index tends to infinity as this index tends to infinity \cite{l-s}. However, the first non-congruence subgroup occurs at index 7.

We now focus on the particular case when $e_2=e_3=g=0$, that is we consider the modular subgroups that are torsion-free and genus zero. They can be generated by parabolic elements only and the Riemann-Hurwitz formula takes the simpler form
\begin{equation}\label{r-h2}
\mu\,=\,6(h-2).
\end{equation}
The torsion-free genus zero congruence subgroups of $\PS$ have been classified into 33 conjugacy classes \cite{seb1}. Beside $\G(2)$ and $\G_0(4)$ which are of index 6, the others are of index a multiple of 12 and the largest index is 60. As we
will see, the smallest index for non-congruence  subgroups in this category is 18. A key tool to analyze the existence of such a modular subgroup is the notion of cusp split, that is of partitioning the index $\mu$ into $h$ cusp widths. To recognize whether a partition corresponds to a modular subgroup is not trivial, but for a given partition, this can be manually  achieved using the permutation groups or graph theory as we will see in later sections.

While the interest in congruence subgroups is self-evident because of their arithmetic properties, there is a growing interest in the non-congruence ones and their modular forms \cite{at-s-d,scholl,long,w-lee}.

\section{Rooted and Unrooted Maps}\label{sec2}
A planar map (or simply a map) is a 2-cell embedding of a connected graph on a 2-dimensional sphere. The elements of the map are the vertices, the edges and the faces which are of dimension 0, 1 and 2 respectively. By a 2-cell embedding we mean that the faces are homeomorphic to discs. Loops and multiple edges are allowed, and  the degree or the valency of a vertex is the number of edges incident with it (loops counting twice). An isomorphism of two maps is a homeomorphism of their spheres such that the image of each element of the first map is an element of the second map with the same dimension. In particular, such isomorphisms preserve the orientation of the maps. An isomorphism of a map onto itself is an automorphism.

Every edge has two ends called darts or semi-edges. They consist of a vertex, a half-edge incident to it and an orientation away from the vertex. Each edge gives rise to two semi-edges (essentially a semi-edge is an edge with an orientation). A rooted map is a map with a choice of one of its semi-edges as a root, and it is considered up to a root-preserving isomorphism. Rooting endows the map with a positive sense of description and right and left sides are specified. Moreover, the enumeration is simpler because the possibility of the map being symmetrical is ruled out. According to Tutte \cite{tut}, a rooted map  has no non-trivial automorphism. An unrooted map is a map  with no assigned root, and they are also considered up to an isomorphism.

\begin{exple}{\rm
	Unrooted maps with 2 edges:
	
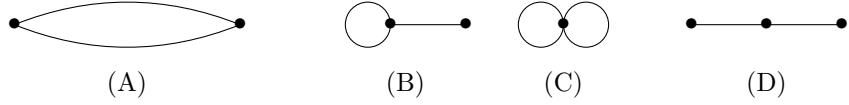
\begin{center}
		\begin{tikzpicture}
		\node  at (0,0) {\textbullet};\node  at (3,0) {\textbullet};
		\draw (0,0) .. controls (1,0.4) and (2,0.4) .. (3,0);
		\draw (0,0) .. controls (1,-0.4) and (2,-0.4) .. (3,0);
		\node at (1.5,-0.8) {(A)};
		\draw (4.7,0) circle (0.3cm); \node  at (5,0) {\textbullet};\draw (5,0) -- (6,0);
		\node  at (6,0) {\textbullet};
		\node at (5.2,-0.8) {(B)};
		\draw (7,0) circle (0.3cm);\node  at (7.3,0) {\textbullet};\draw (7.6,0) circle (0.3cm);
		\node at (7.3,-0.8) {(C)};
		\draw (9,0) -- (11,0);\node  at (9,0) {\textbullet};\node  at (10,0) {\textbullet};
		\node  at (11,0) {\textbullet};
		\node at (10,-0.8) {(D)};
		\end{tikzpicture}
	\end{center}

Notice that by sliding the circle around the sphere, one can see that the map (B) is isomorphic to the map
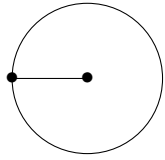
\begin{center}
	\begin{tikzpicture}
		\draw (0,0) circle (1cm); \node  at (-1,0) {\textbullet};\draw (-1,0) -- (0,0);
	\node  at (0,0) {\textbullet};
	\end{tikzpicture}
\end{center}
If we assign a root to these maps, we obtain the following possibilities:

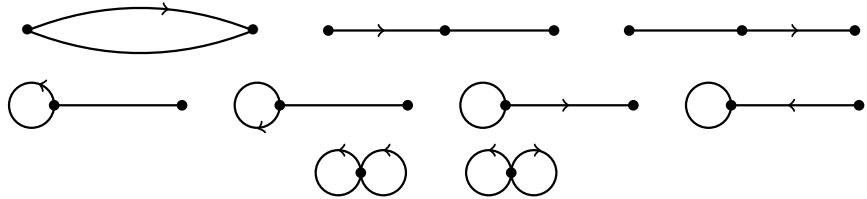
\begin{center}
\begin{tikzpicture}
	\node  at (0,0) {\textbullet};\node  at (3,0) {\textbullet};
\draw[line width=0.3mm,
decoration={markings, mark=at position 0.625 with {\arrow{>}}},
postaction={decorate}
] (0,0) .. controls (1,0.4) and (2,0.4) .. (3,0);
\draw[line width=0.3mm] (0,0) .. controls (1,-0.4) and (2,-0.4) .. (3,0);
%%%%
 \draw[ line width=0.3mm,
 decoration={markings, mark=at position 0.25 with {\arrow{>}}},
postaction={decorate}
]
(4,0) -- (7,0);

% draw the two black dots
\fill (4,0) circle (0.071); 
\fill (5.55,0) circle (0.071); 
\fill (7,0) circle (0.071); 
%%%%
 \draw[ line width=0.3mm,
 decoration={markings, mark=at position 0.75 with {\arrow{>}}},
postaction={decorate}
]
(8,0) -- (11,0);

% draw the two black dots
\fill (8,0) circle (0.071); 
\fill (9.5,0) circle (0.071); 
\fill (11,0) circle (0.071); 
\end{tikzpicture}
\end{center}
\begin{center}
	\begin{tikzpicture}

 \draw[line width=0.3mm,
 decoration={markings, mark=at position 0.2 with {\arrow{>}}},
 postaction={decorate}
 ]
 (0,0) circle (0.3);
 \draw[ line width=0.3mm,
% decoration={markings, mark=at position 0.5 with {\arrow{>}}},
 %postaction={decorate}
 ]
 (.3,0) -- (2,0);
 
 % draw the two black dots
 \fill (0.3,0) circle (0.071); 
 \fill (2,0) circle (0.071); 
 %%%%%%
 \draw[line width=0.3mm,
 decoration={markings, mark=at position 0.8 with {\arrow{<}}},
 postaction={decorate}
 ]
 (3,0) circle (0.3);
 \draw[ line width=0.3mm,
% decoration={markings, mark=at position 0.5 with {\arrow{>}}},
 %postaction={decorate}
 ]
 (3.3,0) -- (5,0);
 
 % draw the two black dots
 \fill (3.3,0) circle (0.071); 
 \fill (5,0) circle (0.071); 
 %%%%%%
 \draw[line width=0.3mm,
 %decoration={markings, mark=at position 0.625 with {\arrow{>}}},
 %postaction={decorate}
 ]
 (6,0) circle (0.3);
 \draw[ line width=0.3mm,
 decoration={markings, mark=at position 0.5 with {\arrow{>}}},
 postaction={decorate}
 ]
 (6.3,0) -- (8,0);
 
 % draw the two black dots
 \fill (6.3,0) circle (0.071); 
 \fill (8,0) circle (0.071); %%%%
  %%%%%%
 \draw[line width=0.3mm,
 %decoration={markings, mark=at position 0.625 with {\arrow{>}}},
 %postaction={decorate}
 ]
 (9,0) circle (0.3);
 \draw[ line width=0.3mm,
 decoration={markings, mark=at position 0.5 with {\arrow{<}}},
 postaction={decorate}
 ]
 (9.3,0) -- (11,0);
 % draw the two black dots
 \fill (9.3,0) circle (0.071); 
 \fill (11,0) circle (0.071); 
 %%%%
 \end{tikzpicture}
\end{center}
\begin{center}
\begin{tikzpicture}
  \draw[line width=0.3mm,
 decoration={markings, mark=at position 0.25 with {\arrow{>}}},
 postaction={decorate}
 ]
 (0,0) circle (0.3);
  \draw[line width=0.3mm,
 decoration={markings, mark=at position 0.25 with {\arrow{>}}},
 postaction={decorate}
 ]
 (0.6,0) circle (0.3);
 \fill (0.3,0) circle (0.071); 
 \draw[line width=0.3mm,
 decoration={markings, mark=at position 0.25 with {\arrow{>}}},
 postaction={decorate}
 ]
 (2,0) circle (0.3);
 \draw[line width=0.3mm,
 decoration={markings, mark=at position 0.25 with {\arrow{<}}},
 postaction={decorate}
 ]
 (2.6,0) circle (0.3);
 \fill (2.3,0) circle (0.071); 
 \end{tikzpicture}
\end{center}
}%%%end rm
\end{exple}
\begin{exple}{\rm If we consider the class of maps with two vertices, each of degree 3, we have:
		
	
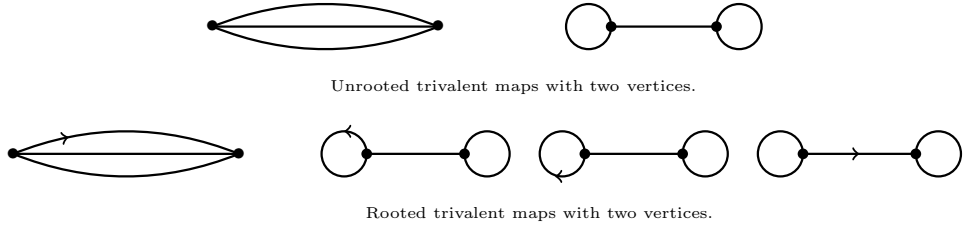
\begin{center}
		\begin{tikzpicture}
		
	\node  at (0,0) {\textbullet};\node  at (3,0) {\textbullet};
	\draw[line width=0.3mm] (0,0) .. controls (1,0.4) and (2,0.4) .. (3,0);
	\draw[line width=0.3mm] (0,0)--(3,0);
	\draw[line width=0.3mm] (0,0) .. controls (1,-0.4) and (2,-0.4) .. (3,0);

%%%%%
\draw[line width=0.3mm,
%decoration={markings, mark=at position 0.625 with {\arrow{>}}},
%postaction={decorate}
]
(5,0) circle (0.3);
\draw[ line width=0.3mm,
%decoration={markings, mark=at position 0.5 with {\arrow{>}}},
%postaction={decorate}
]
(5.3,0) -- (6.7,0);
\draw[line width=0.3mm,
%decoration={markings, mark=at position 0.625 with {\arrow{>}}},
%postaction={decorate}
]
(7,0) circle (0.3);
% draw the two black dots
\fill (5.3,0) circle (0.071); 
\fill (6.7,0) circle (0.071); %%%%
	\node at (4,-0.8) {\tiny Unrooted trivalent maps with two vertices.};
	\end{tikzpicture}
 
\end{center}

 	\begin{center}
 	\begin{tikzpicture}
 	
 	\node  at (0,0) {\textbullet};\node  at (3,0) {\textbullet};
 	\draw[line width=0.3mm,	decoration={markings, mark=at position 0.25 with {\arrow{>}}},
 	postaction={decorate}] (0,0) .. controls (1,0.4) and (2,0.4) .. (3,0);
 	\draw[line width=0.3mm] (0,0)--(3,0);
 	\draw[line width=0.3mm] (0,0) .. controls (1,-0.4) and (2,-0.4) .. (3,0);
 	
 	%%%%%
 	\draw[line width=0.3mm,
 	decoration={markings, mark=at position 0.25 with {\arrow{>}}},
 	postaction={decorate}
 	]
 	(4.4,0) circle (0.3);
 	\draw[ line width=0.3mm,
 	%decoration={markings, mark=at position 0.5 with {\arrow{>}}},
 	%postaction={decorate}
 	]
 	(4.7,0) -- (6,0);
 	\draw[line width=0.3mm,
 	%decoration={markings, mark=at position 0.625 with {\arrow{>}}},
 	%postaction={decorate}
 	]
 	(6.3,0) circle (0.3);
 	% draw the two black dots
 	\fill (4.7,0) circle (0.071); 
 	\fill (6,0) circle (0.071); %%%%
 		%%%%%
 	\draw[line width=0.3mm,
 	decoration={markings, mark=at position 0.75 with {\arrow{<}}},
 	postaction={decorate}
 	]
 	(7.3,0) circle (0.3);
 	\draw[ line width=0.3mm,
 	%decoration={markings, mark=at position 0.5 with {\arrow{>}}},
 	%postaction={decorate}
 	]
 	(7.6,0) -- (8.9,0);
 	\draw[line width=0.3mm,
 	%decoration={markings, mark=at position 0.625 with {\arrow{>}}},
 	%postaction={decorate}
 	]
 	(9.2,0) circle (0.3);
 	% draw the two black dots
 	\fill (7.6,0) circle (0.071); 
 	\fill (8.9,0) circle (0.071); %%%%
 		%%%%%
 	\draw[line width=0.3mm,
 	%decoration={markings, mark=at position 0.25 with {\arrow{>}}},
 	%postaction={decorate}
 	]
 	(10.2,0) circle (0.3);
 	\draw[ line width=0.3mm,
 	decoration={markings, mark=at position 0.5 with {\arrow{>}}},
 	postaction={decorate}
 	]
 	(10.5,0) -- (12,0);
 	\draw[line width=0.3mm,
 	%decoration={markings, mark=at position 0.625 with {\arrow{>}}},
 	%postaction={decorate}
 	]
 	(12.3,0) circle (0.3);
 	% draw the two black dots
 	\fill (10.5,0) circle (0.071); 
 	\fill (12,0) circle (0.071); %%%%
 	\node at (7,-0.8) {\tiny Rooted trivalent maps with two vertices.};
 	\end{tikzpicture}
 	
 \end{center}
}%%%%end rm
\end{exple}

Following the classical notation in enumerating rooted and unrooted maps of a given class, we use $M'$ for the number of rooted maps and $M^+$ for the number of non-isomorphic unrooted maps of the given class.
According to \cite{walsh}, any automorphism of a map is defined by its action on the semi-edges as it preserves the orientation. Moreover, such action
is regular in the sense that as a permutation of the elements of the map, it consists of disjoint cycles of the same length. Suppose $\mathcal{C}(n)$ is a set of  maps with $n$ edges (and $2n$ semi-edges) of a given class. This class can be specified in many ways such as the valency of each vertex or the shape of the faces etc. Moreover, suppose  that $\mathcal{C}(n)$ satisfies the following {\em closure condition}: if the semi-edges of a map $G\in \mathcal{C}(n)$ are labeled, 
then for   $\sigma\in S_{2n}$, $\sigma G$ is also in $\mathcal{C}(n)$. In other words, if $G$ is an unrooted and an unlabeled map, then all maps obtained by  labeling  $G$ simultaneously belong or do not belong to $\mathcal{C}(n)$. Under the closure condition and 
 using the  Burnside Lemma \cite{b1} we have \cite{lisk85}
\begin{equation}\label{burnside}
M^+(n)\,=\, \frac{1}{(2n)!}\,\sum _{\sigma\in S_{2n}}\,m(\sigma,n)
\end{equation}
where $m(\sigma,n)$ is the number of maps in $\mathcal{C}(n)$ whose edges have been labeled and that are invariant under $\sigma$. Since the action of $S_{2n}$ is regular, we have $m(\sigma,n)=0$ if $\sigma$ is not regular. As for the trivial permutation $\epsilon$, we have
\[
m(\epsilon,n)=(2n-1)!M'(n)
\]
because a labeled map is considered as a rooted map with its first semi-edge as its root and there are $(2n-1)!$ ways to label the remaining $2n-1$ semi-edges. Therefore
\[
M^+(n)=\frac{M'(n)}{2n}+\frac{1}{(2n)!}\,\sum\,m(\sigma,n),
\]
where the sum is taken over the non-trivial regular permutations of $S_{2n}$. This  key formula  is used in enumerating unrooted maps as it is simpler  to find closed formulas for $M'(n)$ for a given class and to evaluate $m(\sigma,n)$ using the labeling or the face coloring. As an example , if we consider the class of all maps with $n$ edges, then Tutte  has shown in \cite{tut} that
\[
M'(n)\,=\, \frac{2\cdot 3^n\cdot(2n)!}{n!(n+2)!.}
\]
Furthermore, Liskovets has shown in \cite{lisk85} that 
\begin{align*}
	M^+(n)=&\,\frac{1}{2n}\left[M'(n)+\sum_{t|n,\,t<n}\,\phi(n/t)(t+2)(t+1)M'(t) \right]\\
    	&+\left\{
    	\begin{array}{ll}
		&\frac{n+3}{4}M'\left(\frac{n-1}{2}\right)\quad \mbox{if }n  \mbox{ is odd},\\ \\
		&\frac{n-3}{4}M'\left(\frac{n-2}{2}\right)\quad \mbox{if }n \mbox{ is even},
		\end{array}
	\right.
	\end{align*}
where $\phi$ is the totient Euler function.

We now focus on the class of maps with specified valency at the vertices as their enumeration has been the subject of several articles. Indeed many attempts have been made  to find closed or recursive formulas for the number of  maps with $n$ edges which have $r_i$ vertices of degree $i$, $i\geq1$, where the total degree is $2n=\sum\, ir_i$. A particular case consists of the regular maps having the same degree at all the vertices. For this class, closed recursive formulas have been obtained for Eulerian or even-valent maps \cite{gao,lisk85,lisk,mul,tut,wormald}.  
However, for regular odd-valent maps, the situation is not as simple except for the case $n=3$ where we have a closed formula for enumerating the trivalent rooted maps due to Mullin \cite{mul}, and a formula for enumerating the trivalent unrooted maps due to Liskovets \cite{lisk}. Since the total degree of such maps is even, the number of vertices must be even, say $2p$, with  a positive integer $p$ and the number of edges is then $3p$. From \cite{mul}, we have
\begin{align}\label{rooted-triv}
M'(3p)&=\left\{\begin{array}{ll}
	&\displaystyle \frac{(3p/2)!2^{3p}}{(p+1)!((p+2)/2)!} \quad \mbox{if  } p \mbox{ is even}\\ \\
	&\displaystyle\frac{3(3p+1)!((p+1)/2)!2^p}{((3p+3)/2)!(p+2)!p!}\quad \mbox{if } p \mbox{ is odd}.
	\end{array}
	\right.
\end{align}
As for the unrooted trivalent maps, we have $M^+(3)=2$ and for $p\geq 2$, we have
\begin{align}\label{unroot-triv}
	M^+(3p)=&\frac{1}{6p}\left[M'(3p)+\frac{1}{2}\sum_{t|p,\,t<p} \,\phi(p/t)(t+2)(t+1)M'(3t) \right]\nonumber \\  
	  &+ \frac{2}{3}\left\{
	  \begin{array}{ll}
	  \displaystyle	(p-2)M'(p-4),\  &\mbox{ if } 3|(p-1)\\
	 \displaystyle 	\frac{p+4}{3}M'(p-2),\ &\mbox{ if }3|(p+1)\\
	    	0\,, &\mbox{ if }3|p
	  \end{array}
	  \right.\\ 
	&	+\frac{1}{4}\left\{
		\begin{array}{ll}
	\displaystyle		(p+3)M'\left(\frac{3p-3}{2}\right)\,,&\mbox{ if } p\mbox{ is odd}\\
	\displaystyle	(3p-2)M'\left(\frac{3p-6}{2}\right) \,,&\mbox{ if } p\mbox{ is even}\,.
			\end{array}
		\right.\nonumber 
\end{align}

\section{Enumeration of modular subgroups}

In this section, we will apply the above enumeration schemes to the torsion-free genus zero modular subgroups. Let $\G$ be such a group of index $\mu=6(h-2)$ where $h$ is the parabolic class number. Let $X$ be the set of $\G-$cosets. According to \secref{sec1} we associate to $\G$ a unique 
equivalence class $[\phi,\psi]_{x_1}$ such that the permutations $\phi$ and $\psi$ of $X$ satisfy $\phi^2=\psi^3=1$, $x_1$ is a fixed element of $X$ with $\langle\phi,\psi\rangle$ acting transitively on $X$ and the equivalence class is taken modulo the existence of a permutation $\lambda$ satisfying \eqref{equiv}. To be more precise, consider the images in $\PS$ of the matrices 
$\displaystyle S=\binom{0\ \ -1}{1\ \  \ \ \ 0}$ and $\displaystyle V=\binom{1\  -1}{1\ \ \ \ 0}$. Then $\langle S,V\rangle=\PS$ and $\phi$ and $\psi$ can be taken  respectively as $S$ and $V$ acting on the cosets and $x_1$ can be the coset $\G$. We can  now form the following Schreier graph which is a planar directed graph that summarizes the action of $\phi$ and $\psi$ on $X$: The set of vertices is $X$. If $\phi x=y$ then $x$ and $y$ are joined by a directed blue edge. Since $\phi^2=1$, whenever $x$ and $y$ are connected, then they must be connected with edges with two opposite directions. Similarly, if  $\psi x=y$ then $x$ and $y$ are connected by a directed red edge. Since $\psi^3=1$, each red edge belongs to an oriented triangle and we chose the orientation to be counter-clockwise.
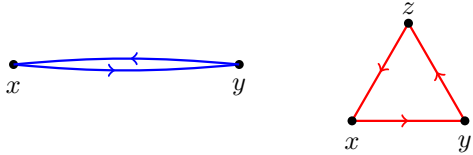
\begin{center}
	\begin{tikzpicture}[scale=1.5]
	
		\draw[line width=0.3mm,	decoration={markings, mark=at position 0.45 with {\arrow{>}}},
		postaction={decorate},color=blue
		] (0,0) .. controls (0.75,-0.07) and (1.25,-0.07) .. (2,0);
		\fill (0,0) circle (0.04); \fill (2,0) circle (0.04); 
		\draw[line width=0.3mm,	decoration={markings, mark=at position 0.55 with {\arrow{<}}},
		postaction={decorate},color=blue
		] (0,0) .. controls (0.75,0.07) and (1.25,0.07) .. (2,0);
	\node at (0,-.2) {$x$};	\node at (2,-.2) {$y$};
	
	\draw[line width=0.3mm,	decoration={markings, mark=at position 0.5 with {\arrow{>}}},
	postaction={decorate},color=red] (3,-0.5)--(4,-0.5);
	\draw[line width=0.3mm,	decoration={markings, mark=at position 0.5 with {\arrow{>}}},
	postaction={decorate},color=red] (4,-0.5)--(3.5,0.366);
	\draw[line width=0.3mm,	decoration={markings, mark=at position 0.5 with {\arrow{>}}},
	postaction={decorate},color=red] (3.5,0.366)-- (3,-0.5);
		\fill (3,-0.5) circle (0.04); \fill (4,-0.5) circle (0.04); 	\fill (3.5,0.366) circle (0.04);
			\node at (3,-.7) {$x$};	\node at (4,-.7) {$y$};	\node at (3.5,0.5) {$z$};	
	\end{tikzpicture}
\end{center}
In other words, the blue arrow indicates the action of $\phi$ and the red one indicates the action of $\psi$. The fact that $\G$ is torsion-free means that the directed graph has no loops, and the genus zero conditions means that the graph is planar (with no edges crossing).

For example,  the level 2 principal congruence group $\G(2)$  is of index 6 and has  the following graph:
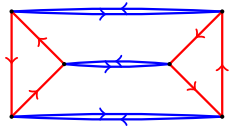
\begin{center}
	\begin{tikzpicture}[scale=0.7
	]
		\draw[line width=0.3mm,	decoration={markings, mark=at position 0.45 with {\arrow{>}}},
	postaction={decorate},color=blue
	] (0,-1) .. controls (1.5,-0.93) and (2.5,-0.93) .. (4,-1);
	\draw[line width=0.3mm,	decoration={markings, mark=at position 0.55 with {\arrow{<}}},
	postaction={decorate},color=blue
	] (0,-1) .. controls (1.5,-1.07) and (2.5,-1.07) .. (4,-1);
	%%%%%
		\draw[line width=0.3mm,	decoration={markings, mark=at position 0.45 with {\arrow{>}}},
	postaction={decorate},color=blue
	] (0,1) .. controls (1.5,0.93) and (2.5,0.93) .. (4,1);
	\draw[line width=0.3mm,	decoration={markings, mark=at position 0.55 with {\arrow{<}}},
	postaction={decorate},color=blue
	] (0,1) .. controls (1.5,1.07) and (2.5,1.07) .. (4,1);
	%%%%%
	\draw[line width=0.3mm,	decoration={markings, mark=at position 0.45 with {\arrow{>}}},
	postaction={decorate},color=blue
	] (1,0) .. controls (1.75,-0.07) and (2.25,-0.07) .. (3,0);
	\draw[line width=0.3mm,	decoration={markings, mark=at position 0.55 with {\arrow{<}}},
	postaction={decorate},color=blue
	] (1,0) .. controls (1.75,0.07) and (2.25,0.07) .. (3,0);
%%%%%%%
	\draw[line width=0.3mm,	decoration={markings, mark=at position 0.5 with {\arrow{>}}},
	postaction={decorate},color=red] (0,1)--(0,-1);
	\draw[line width=0.3mm,	decoration={markings, mark=at position 0.5 with {\arrow{>}}},
	postaction={decorate},color=red] (0,-1)--(1,0);
	\draw[line width=0.3mm,	decoration={markings, mark=at position 0.5 with {\arrow{>}}},
	postaction={decorate},color=red] (1,0)-- (0,1);
	\fill (0,1) circle (0.04); \fill (0,-1) circle (0.04); 	\fill (1,0) circle (0.04);
	%%%
	\draw[line width=0.3mm,	decoration={markings, mark=at position 0.5 with {\arrow{>}}},
postaction={decorate},color=red] (3,0)--(4,-1);
\draw[line width=0.3mm,	decoration={markings, mark=at position 0.5 with {\arrow{>}}},
postaction={decorate},color=red] (4,-1)--(4,1);
\draw[line width=0.3mm,	decoration={markings, mark=at position 0.5 with {\arrow{>}}},
postaction={decorate},color=red] (4,1)-- (3,0);
\fill (4,1) circle (0.04); \fill (4,-1) circle (0.04); 	\fill (3,0) circle (0.04);
	\end{tikzpicture}
\end{center}
Here the permutation $\phi\psi$ decomposes into 3 cycles each of length 2, which corresponds to a cusp split 2-2-2. Another congruence group of index 6 is $\G_0(4)$ and has the following Schreier graph:\

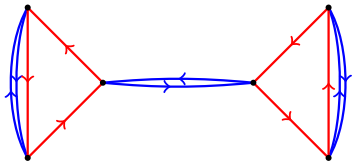
\begin{center}
	\begin{tikzpicture}[scale=1
	]
	\draw[line width=0.3mm,	decoration={markings, mark=at position 0.45 with {\arrow{>}}},
	postaction={decorate},color=blue
	] (0,-1) .. controls (-.3,-0.5) and (-.3,0.5) .. (0,1);
	\draw[line width=0.3mm,	decoration={markings, mark=at position 0.55 with {\arrow{<}}},
	postaction={decorate},color=blue
	] (0,-1) .. controls (-0.2,-0.5) and (-.2,0.5) .. (0,1);
	%%%%%
		\draw[line width=0.3mm,	decoration={markings, mark=at position 0.45 with {\arrow{>}}},
	postaction={decorate},color=blue
	] (4,-1) .. controls (4.2,-0.5) and (4.2,0.5) .. (4,1);
	\draw[line width=0.3mm,	decoration={markings, mark=at position 0.55 with {\arrow{<}}},
	postaction={decorate},color=blue
	] (4,-1) .. controls (4.3,-0.5) and (4.3,0.5) .. (4,1);
	%%%%%
	\draw[line width=0.3mm,	decoration={markings, mark=at position 0.45 with {\arrow{>}}},
	postaction={decorate},color=blue
	] (1,0) .. controls (1.75,-0.07) and (2.25,-0.07) .. (3,0);
	\draw[line width=0.3mm,	decoration={markings, mark=at position 0.55 with {\arrow{<}}},
	postaction={decorate},color=blue
	] (1,0) .. controls (1.75,0.07) and (2.25,0.07) .. (3,0);
	%%%%%%%
	\draw[line width=0.3mm,	decoration={markings, mark=at position 0.5 with {\arrow{>}}},
	postaction={decorate},color=red] (0,1)--(0,-1);
	\draw[line width=0.3mm,	decoration={markings, mark=at position 0.5 with {\arrow{>}}},
	postaction={decorate},color=red] (0,-1)--(1,0);
	\draw[line width=0.3mm,	decoration={markings, mark=at position 0.5 with {\arrow{>}}},
	postaction={decorate},color=red] (1,0)-- (0,1);
	\fill (0,1) circle (0.04); \fill (0,-1) circle (0.04); 	\fill (1,0) circle (0.04);
	%%%
	\draw[line width=0.3mm,	decoration={markings, mark=at position 0.5 with {\arrow{>}}},
	postaction={decorate},color=red] (3,0)--(4,-1);
	\draw[line width=0.3mm,	decoration={markings, mark=at position 0.5 with {\arrow{>}}},
	postaction={decorate},color=red] (4,-1)--(4,1);
	\draw[line width=0.3mm,	decoration={markings, mark=at position 0.5 with {\arrow{>}}},
	postaction={decorate},color=red] (4,1)-- (3,0);
	\fill (4,1) circle (0.04); \fill (4,-1) circle (0.04); 	\fill (3,0) circle (0.04);
	\end{tikzpicture}
\end{center}
Here the cusp split is 4-1-1. In general each vertex is incident to two blue edges and two red edges. It is straightforward that  such a graph characterizes completely and uniquely the class $[\phi,\psi]_{x_1}$. Indeed, one only needs to assign a vertex as $x_1$ and label the others as $x_2,\cdots,x_{\mu}$ and then write down the permutations. Any other labeling would require the action of a permutation $\lambda$ that fixes $x_1$ and conjugates the permutations making sure it will yield the same class $[\phi,\psi]_{x_1}$.

To each Schreier graph we can associate a unique rooted planar map in the following way: The $\psi-$orbits which are the red triangles are the vertices, and there is an edge between two such vertices if there is a blue (double) edge connecting the two orbits. In other words, we shrink every red triangle to a   vertex and the two blue edges become one undirected edge. In the meantime, the root vertex is the orbit of $x_1$ and the oriented edge is simply the one coming from the blue edge incident to $x_1$ with the orientation away from $x_1$.

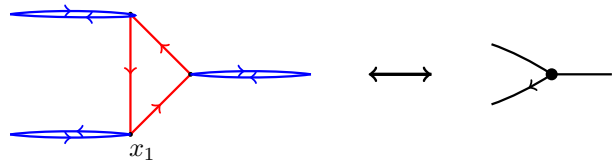
\begin{center}
	\begin{tikzpicture}[scale=0.8
	]
	\draw[line width=0.3mm,	decoration={markings, mark=at position 0.5 with {\arrow{>}}},
postaction={decorate},color=red] (0,1)--(0,-1);
\draw[line width=0.3mm,	decoration={markings, mark=at position 0.5 with {\arrow{>}}},
postaction={decorate},color=red] (0,-1)--(1,0);
\draw[line width=0.3mm,	decoration={markings, mark=at position 0.5 with {\arrow{>}}},
postaction={decorate},color=red] (1,0)-- (0,1);
\fill (0,1) circle (0.04); \fill (0,-1) circle (0.04); 	\fill (1,0) circle (0.04);
%%%%
	\draw[line width=0.3mm,	decoration={markings, mark=at position 0.45 with {\arrow{>}}},
postaction={decorate},color=blue
] (0,1) .. controls (00.5,0.93) and (-1.5,0.93) .. (-2,1);
\draw[line width=0.3mm,	decoration={markings, mark=at position 0.55 with {\arrow{<}}},
postaction={decorate},color=blue
] (0,1) .. controls (-0.5,1.07) and (-1.5,1.07) .. (-2,1);
%%%%%%%
\draw[line width=0.3mm,	decoration={markings, mark=at position 0.45 with {\arrow{>}}},
postaction={decorate},color=blue
] (0,-1) .. controls (-0.5,-0.93) and (-1.5,-0.93) .. (-2,-1);
\draw[line width=0.3mm,	decoration={markings, mark=at position 0.55 with {\arrow{<}}},
postaction={decorate},color=blue
] (0,-1) .. controls (-0.5,-1.07) and (-1.5,-1.07) .. (-2,-1);
%%%%%%%
\draw[line width=0.3mm,	decoration={markings, mark=at position 0.45 with {\arrow{>}}},
postaction={decorate},color=blue
] (1,0) .. controls (1.5,0.07) and (2.5,0.07) .. (3,0);
\draw[line width=0.3mm,	decoration={markings, mark=at position 0.55 with {\arrow{<}}},
postaction={decorate},color=blue
] (1,0) .. controls (1.5,-0.07) and (2.5,-0.07) .. (3,0);
\node at (0.2,-1.3) {$x_1$};
%%%
\draw[line width=0.4mm,	decoration={markings, mark=at position 1 with {\arrow{>}}},
postaction={decorate}] (4,0)--(5,0);
\draw[line width=0.4mm,	decoration={markings, mark=at position 0.1 with {\arrow{<}}},
postaction={decorate}] (4,0)--(5,0);
%%%
	\draw[line width=0.3mm,	decoration={markings, mark=at position 0.4 with {\arrow{>}}},
postaction={decorate}
] (7,0) .. controls (6.7,-0.2) and (6.3,-0.4) .. (6,-0.5); \fill (7,0) circle (0.1);
	\draw[line width=0.3mm
] (7,0) .. controls (6.7,0.2) and (6.3,0.4) .. (6,0.5); 
\draw[line width=0.3mm] (7,0)--(8,0);
\end{tikzpicture}
\end{center}
It is clear that the resulting map is trivalent.
The process is easily reversible because once we change the root vertex into a triangle, then $x_1$ is positioned at the vertex of the triangle that is incident to the double edge that replaced the root edge.  As examples, the above Schreier graphs with a choice of the position of $x_1$ would correspond to:

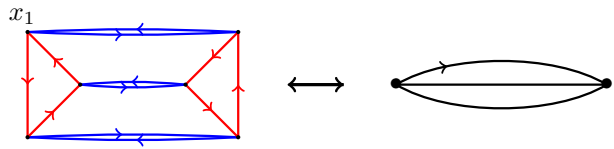
\begin{center}
	\begin{tikzpicture}[scale=0.7
	]
	\node at (-0.1,1.3) {$x_1$};
	
	\draw[line width=0.3mm,	decoration={markings, mark=at position 0.45 with {\arrow{>}}},
	postaction={decorate},color=blue
	] (0,-1) .. controls (1.5,-0.93) and (2.5,-0.93) .. (4,-1);
	\draw[line width=0.3mm,	decoration={markings, mark=at position 0.55 with {\arrow{<}}},
	postaction={decorate},color=blue
	] (0,-1) .. controls (1.5,-1.07) and (2.5,-1.07) .. (4,-1);
	%%%%%
	\draw[line width=0.3mm,	decoration={markings, mark=at position 0.45 with {\arrow{>}}},
	postaction={decorate},color=blue
	] (0,1) .. controls (1.5,0.93) and (2.5,0.93) .. (4,1);
	\draw[line width=0.3mm,	decoration={markings, mark=at position 0.55 with {\arrow{<}}},
	postaction={decorate},color=blue
	] (0,1) .. controls (1.5,1.07) and (2.5,1.07) .. (4,1);
	%%%%%
	\draw[line width=0.3mm,	decoration={markings, mark=at position 0.45 with {\arrow{>}}},
	postaction={decorate},color=blue
	] (1,0) .. controls (1.75,-0.07) and (2.25,-0.07) .. (3,0);
	\draw[line width=0.3mm,	decoration={markings, mark=at position 0.55 with {\arrow{<}}},
	postaction={decorate},color=blue
	] (1,0) .. controls (1.75,0.07) and (2.25,0.07) .. (3,0);
	%%%%%%%
	\draw[line width=0.3mm,	decoration={markings, mark=at position 0.5 with {\arrow{>}}},
	postaction={decorate},color=red] (0,1)--(0,-1);
	\draw[line width=0.3mm,	decoration={markings, mark=at position 0.5 with {\arrow{>}}},
	postaction={decorate},color=red] (0,-1)--(1,0);
	\draw[line width=0.3mm,	decoration={markings, mark=at position 0.5 with {\arrow{>}}},
	postaction={decorate},color=red] (1,0)-- (0,1);
	\fill (0,1) circle (0.04); \fill (0,-1) circle (0.04); 	\fill (1,0) circle (0.04);
	%%%
	\draw[line width=0.3mm,	decoration={markings, mark=at position 0.5 with {\arrow{>}}},
	postaction={decorate},color=red] (3,0)--(4,-1);
	\draw[line width=0.3mm,	decoration={markings, mark=at position 0.5 with {\arrow{>}}},
	postaction={decorate},color=red] (4,-1)--(4,1);
	\draw[line width=0.3mm,	decoration={markings, mark=at position 0.5 with {\arrow{>}}},
	postaction={decorate},color=red] (4,1)-- (3,0);
	\fill (4,1) circle (0.04); \fill (4,-1) circle (0.04); 	\fill (3,0) circle (0.04);
	%%%%%%
	\draw[line width=0.4mm,	decoration={markings, mark=at position 1 with {\arrow{>}}},
	postaction={decorate}] (5,0)--(6,0);
	\draw[line width=0.4mm,	decoration={markings, mark=at position 0.1 with {\arrow{<}}},
	postaction={decorate}] (5,0)--(6,0);
	%%%%%
		\node  at (7,0) {\textbullet};\node  at (11,0) {\textbullet};
	\draw[line width=0.3mm,	decoration={markings, mark=at position 0.25 with {\arrow{>}}},
	postaction={decorate}] (7,0) .. controls (8,0.6) and (10,0.6) .. (11,0);
	\draw[line width=0.3mm] (7,0)--(11,0);
	\draw[line width=0.3mm] (7,0) .. controls (8,-0.6) and (10,-0.6) .. (11,0);
	
	\end{tikzpicture}
\end{center}
In fact if we change the position of $x_1$, we get the same directed graphs on the sphere % (the triangles are not chiral when facing each other) 
which yields a unique  trivalent rooted map with two vertices.
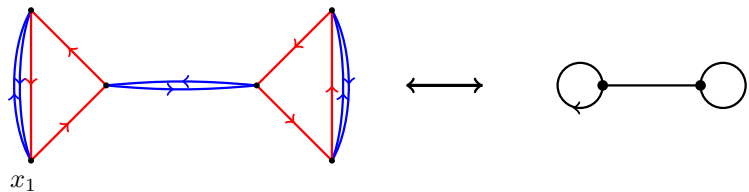
\begin{center}
	\begin{tikzpicture}[scale=1]
	\node at (-0.1,-1.3) {$x_1$};
	
	\draw[line width=0.3mm,	decoration={markings, mark=at position 0.45 with {\arrow{>}}},
	postaction={decorate},color=blue
	] (0,-1) .. controls (-.3,-0.5) and (-.3,0.5) .. (0,1);
	\draw[line width=0.3mm,	decoration={markings, mark=at position 0.55 with {\arrow{<}}},
	postaction={decorate},color=blue
	] (0,-1) .. controls (-0.2,-0.5) and (-.2,0.5) .. (0,1);
	%%%%%
	\draw[line width=0.3mm,	decoration={markings, mark=at position 0.45 with {\arrow{>}}},
	postaction={decorate},color=blue
	] (4,-1) .. controls (4.2,-0.5) and (4.2,0.5) .. (4,1);
	\draw[line width=0.3mm,	decoration={markings, mark=at position 0.55 with {\arrow{<}}},
	postaction={decorate},color=blue
	] (4,-1) .. controls (4.3,-0.5) and (4.3,0.5) .. (4,1);
	%%%%%
	\draw[line width=0.3mm,	decoration={markings, mark=at position 0.45 with {\arrow{>}}},
	postaction={decorate},color=blue
	] (1,0) .. controls (1.75,-0.07) and (2.25,-0.07) .. (3,0);
	\draw[line width=0.3mm,	decoration={markings, mark=at position 0.55 with {\arrow{<}}},
	postaction={decorate},color=blue
	] (1,0) .. controls (1.75,0.07) and (2.25,0.07) .. (3,0);
	%%%%%%%
	\draw[line width=0.3mm,	decoration={markings, mark=at position 0.5 with {\arrow{>}}},
	postaction={decorate},color=red] (0,1)--(0,-1);
	\draw[line width=0.3mm,	decoration={markings, mark=at position 0.5 with {\arrow{>}}},
	postaction={decorate},color=red] (0,-1)--(1,0);
	\draw[line width=0.3mm,	decoration={markings, mark=at position 0.5 with {\arrow{>}}},
	postaction={decorate},color=red] (1,0)-- (0,1);
	\fill (0,1) circle (0.04); \fill (0,-1) circle (0.04); 	\fill (1,0) circle (0.04);
	%%%
	\draw[line width=0.3mm,	decoration={markings, mark=at position 0.5 with {\arrow{>}}},
	postaction={decorate},color=red] (3,0)--(4,-1);
	\draw[line width=0.3mm,	decoration={markings, mark=at position 0.5 with {\arrow{>}}},
	postaction={decorate},color=red] (4,-1)--(4,1);
	\draw[line width=0.3mm,	decoration={markings, mark=at position 0.5 with {\arrow{>}}},
	postaction={decorate},color=red] (4,1)-- (3,0);
	\fill (4,1) circle (0.04); \fill (4,-1) circle (0.04); 	\fill (3,0) circle (0.04);
	%%%%%%
	\draw[line width=0.4mm,	decoration={markings, mark=at position 1 with {\arrow{>}}},
	postaction={decorate}] (5,0)--(6,0);
	\draw[line width=0.4mm,	decoration={markings, mark=at position 0.1 with {\arrow{<}}},
	postaction={decorate}] (5,0)--(6,0);
	%%%%%
	
	%%%%%
	\draw[line width=0.3mm,
	decoration={markings, mark=at position 0.75 with {\arrow{<}}},
	postaction={decorate}
	]
	(7.3,0) circle (0.3);
	\draw[ line width=0.3mm,
	%decoration={markings, mark=at position 0.5 with {\arrow{>}}},
	%postaction={decorate}
	]
	(7.6,0) -- (8.9,0);
	\draw[line width=0.3mm,
	%decoration={markings, mark=at position 0.625 with {\arrow{>}}},
	%postaction={decorate}
	]
	(9.2,0) circle (0.3);
	% draw the two black dots
	\fill (7.6,0) circle (0.071); 
	\fill (8.9,0) circle (0.071); %%%%
	%%%%%
	\end{tikzpicture}
\end{center}
In this case, moving $x_1$ on its $\psi-$orbit gives three different directed graphs to which correspond three different trivalent rooted maps with two vertices that were listed in the previous section. 

Recall that a torsion-free genus zero modular subgroup has index $\mu=6(h-2)$, which is the  number of vertices in the corresponding Schreier graph and so the associated trivalent rooted map has $2(h-2)$ vertices and $3(h-2)$ edges. Moreover, the parabolic class number $h$ corresponds to the number of faces and the cusp split provides the degrees of all the faces, that is, the number of edges incident to each face.

We have so far established a one to one correspondence between the conjugacy classes $[\phi,\psi]_{x_1}$ and the trivalent rooted maps with $2(h-2)$ vertices. Moreover, when we ignore the element $x_1$,  we obtain a correspondence between the classes $[\phi,\psi]$ (or the Schreier graphs with no label $x_1$ ) and the non-isomorphic unrooted trivalent maps with $2(h-2)$ vertices. Notice that $\mu$ is the total degree of the map. We have therefore proved the following:
\begin{thm} If $6 |\mu$, then 
	\begin{enumerate}
		\item there is a one-to-one correspondence between the torsion-free genus zero modular subgroups of index $\mu$ and the trivalent rooted maps with
		$\mu/3$ vertices and
		\item there is a one-to-one correspondence between the conjugacy classes of  torsion-free genus zero modular subgroups of index $\mu$ and the trivalent unrooted maps with  $\mu/3$ vertices.
	\end{enumerate}
\end{thm}
As a consequence and using the previous section, we have:
\begin{thm} If $6 | \mu$, then the number  of  torsion-free genus zero modular subgroups of index $\mu$ is given by:

\begin{align*}
	&\displaystyle N(\mu)=\frac{(\mu/4)!2^{\mu/2}}{(\frac{\mu}{6}+1)!(\frac{\mu}{12}+1)!} \quad \mbox{if  } 12 | \mu\,,\\ 
	&\displaystyle N(\mu)=\frac{3(\frac{\mu}{2}+1)!(\frac{\mu}{12}+\frac{1}{2})!2^{\mu/6}}{\left(\frac{\mu}{4}+\frac{3}{2}\right)!(\frac{\mu}{6}+2)!(\frac{\mu}{6})!}\quad \mbox{if  } 12 \nmid\mu.
	\end{align*}
If $6\nmid \mu$, then $N(\mu)=0$.
\end{thm}

As for the conjugacy classes of these groups, we have:
\begin{thm}
	 If $6 | \mu$, then the number  of  conjugacy classes of torsion-free genus zero modular subgroups of index $\mu$ is given by $N^+(6)=2$ and for $\mu>6$:
	 \begin{align*}\label{unroot-triv}
	 	N^+(\mu)=&\frac{1}{\mu}\left[N(\mu)+\frac{1}{2}\sum_{t|\frac{\mu}{6},\, t<\frac{\mu}{6}} \,\phi(\frac{\mu}{6t})(t+2)(t+1)N(6t) \right] \\  
	 	&+ \frac{2}{3}\left\{
	 	\begin{array}{ll}
	 	\displaystyle	(\frac{\mu}{6}-2)N(\frac{\mu}{3}-8),\  &\mbox{ if } \mu\equiv 6\mod 18\\ \\
	 	\displaystyle	(\frac{\mu}{18}+\frac{4}{3})N(\frac{\mu}{3}-4),\ &\mbox{ if  } \mu\equiv 12\mod 18\\ \\
	 		0\,, &\mbox{ if  }\mu\equiv 0\mod 18\\
	 	\end{array}
	 	\right.\\  
	 	&	+\frac{1}{4}\left\{
	 	\begin{array}{ll}
	 \displaystyle		(\frac{\mu}{6}+3)N\left(\frac{\mu}{2}-3\right)\,,&\mbox{ if  } 4\nmid \mu\\ \\
	\displaystyle 		(\frac{\mu}{2}-2)N\left(\frac{\mu}{2}-6\right) \,,&\mbox{ if  } 4\mid \mu\,,
	 	\end{array}
	 	\right.
	 \end{align*}
 with the convention that $N(0)=1$ and $N(k)=0$ if $k<0$.
\end{thm}
As an example, at index 6, we have the  subgroups  $\G(2)$ which is normal in $\PS$, and $\G_0(4)$ and 2 of its conjugates which account for the four rooted trivalent maps with 2 vertices seen earlier. These conjugates are:
\[
	\binom{0\ \  -1}{1\ \ \  \ \ \ 0}\G_0(4)	\binom{0\ \  -1}{1\ \ \ \ 0}^{-1}=\G^0(4)=\left\{\mat\in\SL \ \mid \ b\equiv 0 \mod 4\right\},
	\]
	\[ 
	\binom{1\  \ -1 }{1\ \ \ \ \  0}\G_0(4)\binom{1\  \ -1 }{1\ \ \ \ \  0}^{-1}=\left\{\mat\in\SL \ \mid\  a+b-c-d\equiv 0 \mod 4\right\}.
		\]

		The following table indicates the values of $N(\mu)$ and $N^+(\mu)$ for $\mu\leq 60$.
		
		\begin{center}
		\begin{tabular}{l|l|l|}
				$\mu$&$N(\mu)$&$N^+(\mu)$\\
				\hline
				6&4&2\\
				12&32&6\\
				18&336&26\\
				24&4096&191\\
				30&54912&1904\\
				36&786432&22078\\
				42&11824384&282388\\
				48&184549376 &3848001 \\
				54&2966845440 &54953996 \\
				60& 48855252992& 814302292\\
				\hline
			\end{tabular}
		\end{center}

\begin{remark}{\em
	The number of faces in the map corresponds to the parabolic class number of the corresponding modular subgroup, and the cusp width is exactly the degree of the face, i.e. the number of edges bordering the face. In other words, it easy to write down the cusp split from the map.}
	\end{remark}

\begin{remark}{\em From \cite{seb1} we have the following information:
At index 6, we have two classes of torsion-free genus 0 modular subgroups, namely $\G(2)$ and $\G_0(4)$. At index 12, we have six such groups and all of them are congruence subgroups, namely $\G(3)$, $\G_0(4)\cap\G(2)$, $\G_1(5)$, $\G_0(6)$, $\G_0(8)$ and $\G_0(9)$. As for the groups of index 18, none of them is  congruence, and among the 191 groups of index 24 only 9 are congruence. The last index at which there are congruence subgroups is 60, with the two groups being $\G(5)$ and $\G_0(25)\cap\G_1(5)$. }
\end{remark}

\section{Elliptic surfaces}
We recall that an elliptic surface $X$ over a curve $C$ is a holomorphic map $\phi:X\longrightarrow C$ between a projective surface $X$ over $\BC$ and a smooth curve $C$ such that the fibers of $\phi$ outside a finite number of points of $C$ are smooth connected curves of genus 1. Moreover, we assume the existence of a distinguished global section  so that the general fibers become elliptic curves. The singular fibers are those which are not elliptic curves.  In \cite{kod1,kod2}, Kodaira has classified the fiber types into two infinite families labeled $I_n$ (a cycle of $n$ irreducible curves) and $I_n^*$ (a cycle of $n+5$ irreducible curves) for $n\geq 0$, and six exceptional types. An elliptic surface is called semi-stable if all the singular fibers are of type $I_n$. Furthermore, If $\G'$ is a subgroup of $\SL$ and $\G$ is its image in $\PS$, following Kodaira in \cite{kod2}, one can associate to $\G'$ an elliptic surface over $X(\G)$  whose homological and functional invariants are constructed from the covering $J: X(\G)\longrightarrow X(1)$. In case $\G'$ is torsion-free, the construction of the elliptic surface is carried out in the following explicit way: We start with the product $\fH \times\BC$ and take the quotient by the semi-direct product of $\G'$ and $\BZ^2$ acting as follows: If $\gamma=\mat\in\G'$, $(m,n)\in\BZ^2$  and $(\tau,z)\in\fH\times\BC$ then
\[
(\gamma,m,n)(\tau,z)=\left(\gamma\cdot\tau\ ,\  \frac{z+m\tau+n}{c\tau+d}\right).
\]
This surface extends to an elliptic surface over $X(\G)$ such that the singular fibers are above the cusps. A cusp of the first kind (fixed by a matrix of trace 2) gives rise to a singular fiber of type $I_n$ where $n$ is the cusp width, and a cusp of the second kind (fixed by a matrix of trace -2) gives rise to a fiber of type $I_n^*$. This construction was extended by Shioda  to groups having elliptic elements \cite{shio1,shio2} and to subgroups of $\PS$ by Nori \cite{nori}. These surfaces are referred to as elliptic  modular surfaces. They have the property that their Mordell-Weil groups are finite. If the base curve $C$ has genus 0, the geometric genus is given by $p_g=h/2\ -\ 2$,   and the arithmetic genus $p_a$ of the surface is equal to $p_g$. 

As an example, the index 6 subgroups give rise to three rational elliptic surfaces over the genus zero sphere with three singular fibers of type $[I_2^*,\,I_2,\,I_2]$,
$[I_4^*,I_1,I_1]$ and $[I_4,I_1^*,I_1]$. Meanwhile, at index 12, we have the Beauville surfaces \cite{beau} which are also rational and semi-stable, namely $[I_9,\,I_1,\,I_1,\,I_1]$, $[I_8,\,I_2,\,I_1,\,I_1]$, $[I_6,\,I_3,\,I_2,\,I_1]$,  $[I_5,\,I_5,\,I_1,\,I_1]$,  $[I_4,I_4,I_2,I_2]$ and
$[I_3,\,I_3,\,I_3,\,I_3]$.

 Notice that the Euler characteristic of the surface is the sum of the Euler characteristic of the singular fibers $e(C_v)$,  and these are equal to $m_v$, the number of irreducible components of $C_v$ if it is of type $I_n$, and are equal to $m_v+1$ for all the other types. In particular, $e(I_n)=n$ and $e(I_n^*)=n+6$. Thus the above rational surfaces at index 6 and 12 have Euler characteristic equal to 12.
 
 At index 24, we encounter elliptic surfaces with six semi-stable singular fibers that  have been studied by Miranda and Pearson \cite{m-p}, Shimada \cite{shim}, Shioda \cite{shio3}, and Beukers and Vidunas \cite{b-v}. These are semi-stable elliptic K3 surfaces over ${\mathbb P}^1$ which are extremal in the sense that they have a maximal Picard number which is 20. Indeed, the latter  is given by \cite{kod1}
 \[
 p\,=\, r+2+\sum_{v=1}^h\, (m_v-1).
 \]
 Here $r$ is the rank of the Mordell-Weill group which is always 0 in our case, and hence $p=20$ for the surfaces under consideration .
 For the rest of this section, we focus on  K3 elliptic surfaces over ${\mathbb P}^1$ that have four singular fibers of type $I_n$ and one of type $I_n^*$.  According to \cite{nori}, these are also modular elliptic surfaces and so the attached modular  subgroup is  torsion-free, genus zero and has index 18 (which is the Euler characteristic of the surface minus 6). They are also extremal as the Picard number is 20. Each subgroup correspond to more than one non-isomorphic elliptic surface by moving the $({ }^*)$ along the various singular fibers accounting for the various lift of $\G$ to $\SL$.

There are 26 such groups according to the previous section. The associated  unrooted graphs are given in the table below.  The first example corresponding to the cusp split 14-1-1-1-1 has been studied by Shioda in \cite{shio3} as one of two elliptic surfaces with a maximal fiber type, the other one being the semi-stable case corresponding to the cusp split 19-1-1-1-1-1 at index 24. It is worthwhile to mention that if a K3 surface is elliptically fibered, then the base curve is necessarily of genus 0 \cite{sch-sh}

We notice that there are two non-isomorphic groups having the same cusp split, namely 7-7-2-1-1. Their graphs are chiral and are the only ones which do not possess any sort of symmetry.

%%%% N1
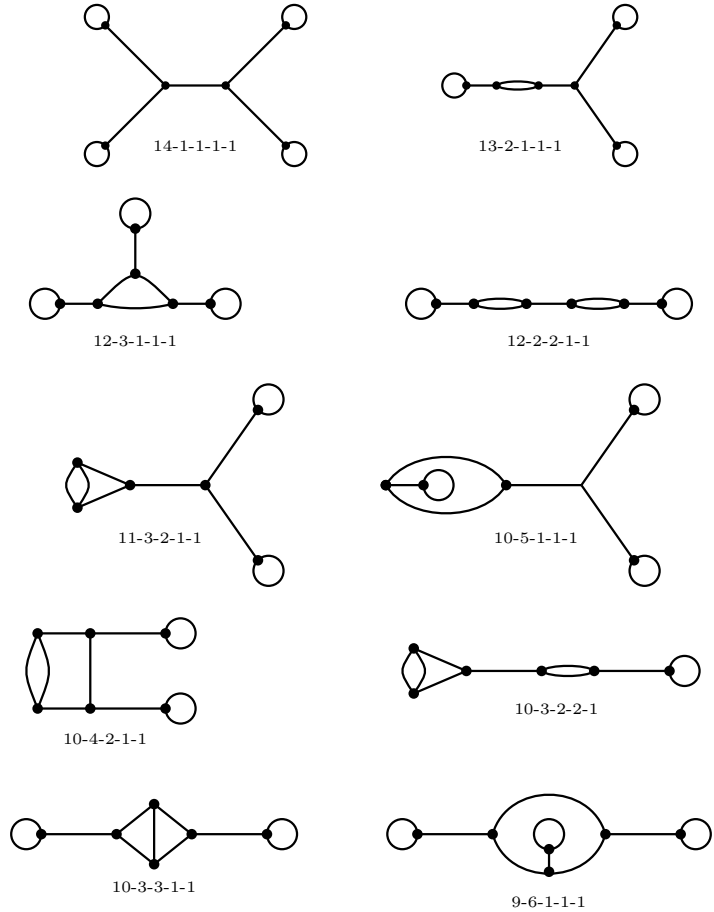
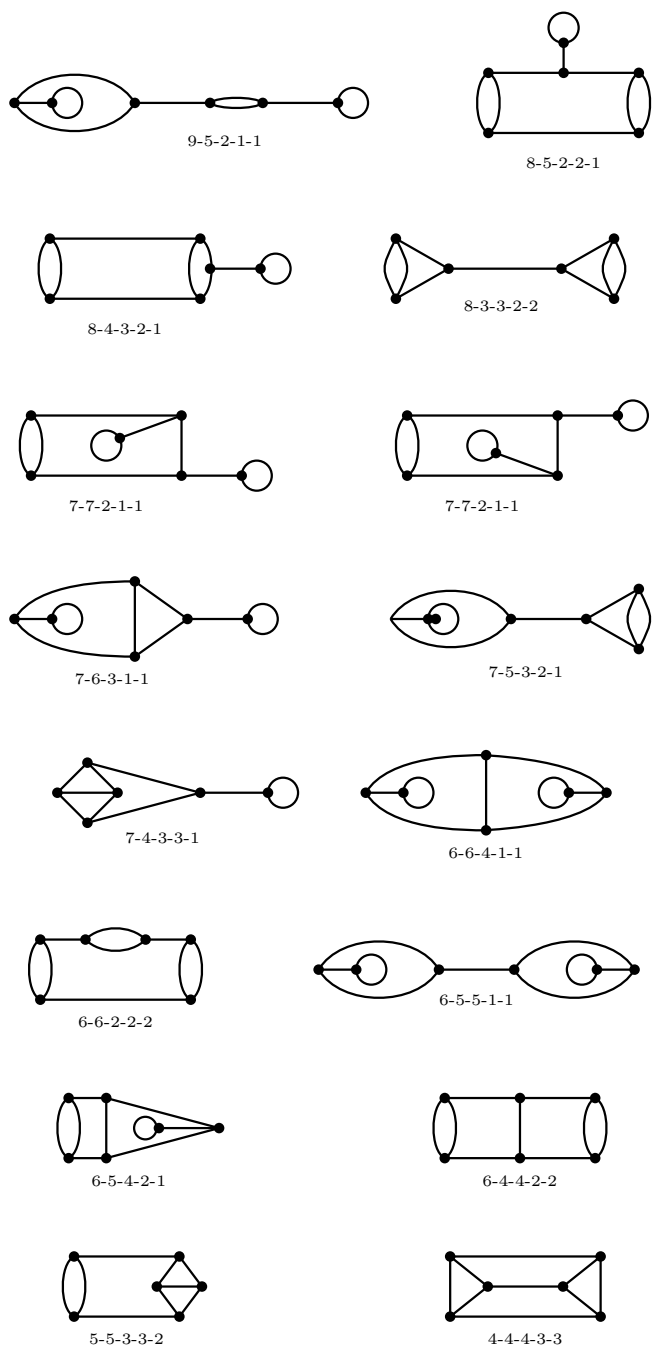
\begin{center}
	\begin{tikzpicture}[scale=0.8]
\draw[line width=0.3mm] (0,0)--(1,0);\draw[line width=0.3mm] (-1,1)--(0,0);\draw[line width=0.3mm] (0,0)--(-1,-1);\draw[line width=0.3mm] (2,1)--(1,0);
\draw[line width=0.3mm] (2,-1)--(1,0);
\fill(0,0) circle (0.07);\fill(-1,1) circle (0.07);\fill(-1,-1) circle (0.07);\fill(1,0) circle (0.07);\fill(2,1) circle (0.07);\fill(2,-1) circle (0.07);
\draw[line width=0.3mm] (-1.14,1.14) circle (0.2);\draw[line width=0.3mm] (2.14,1.14) circle (0.2);
\draw[line width=0.3mm] (-1.14,-1.14) circle (0.2);\draw[line width=0.3mm] (2.14,-1.14) circle (0.2);
{\tiny \node at (0.5,-1){14-1-1-1-1};}
%%%%%%%%%%%%%%%%%%%%%%%%%%
\draw[line width=0.3mm] (4.8,0) circle (0.2);\draw[line width=0.3mm] (5,0)--(5.5,0);
\draw[line width=0.3mm	] (5.5,0) .. controls (5.6,-0.09) and (6.1,-0.09) .. (6.2,0);
\draw[line width=0.3mm	] (5.5,0) .. controls (5.6,0.09) and (6.1,0.09) .. (6.2,0);\draw[line width=0.3mm] (6.2,0)--(6.8,0);
\draw[line width=0.3mm](6.8,0)--(7.5,1);\draw[line width=0.3mm](6.8,0)--(7.5,-1);
\draw[line width=0.3mm] (7.64,1.14) circle (0.2);
\draw[line width=0.3mm] (7.64,-1.14) circle (0.2);
\fill(5,0) circle (0.07);\fill(5.5,0) circle (0.07);\fill(6.2,0) circle (0.07);\fill(6.8,0) circle (0.07);\fill(7.5,1) circle (0.07);\fill(7.5,-1) circle (0.07);
{\tiny \node at (5.9,-1){13-2-1-1-1};}
	\end{tikzpicture}
\end{center}
%%%% N2
\begin{comment}
\begin{center}
	\begin{tikzpicture}[scale=0.8]
	\draw[line width=0.3mm] (-0.2,0) circle (0.2);\draw[line width=0.3mm] (0,0)--(0.5,0);
		\draw[line width=0.3mm	] (0.5,0) .. controls (0.6,-0.09) and (1.1,-0.09) .. (1.2,0);
		\draw[line width=0.3mm	] (0.5,0) .. controls (0.6,0.09) and (1.1,0.09) .. (1.2,0);\draw[line width=0.3mm] (1.2,0)--(1.8,0);
		\draw[line width=0.3mm](1.8,0)--(2.5,1);\draw[line width=0.3mm](1.8,0)--(2.5,-1);
		\draw[line width=0.3mm] (2.64,1.14) circle (0.2);
		\draw[line width=0.3mm] (2.64,-1.14) circle (0.2);
		\fill(0,0) circle (0.07);\fill(0.5,0) circle (0.07);\fill(1.2,0) circle (0.07);\fill(1.8,0) circle (0.07);\fill(2.5,1) circle (0.07);\fill(2.5,-1) circle (0.07);
		{\tiny \node at (0.9,-1){13-2-1-1-1};}
	\end{tikzpicture}
\end{center}
\end{comment}

%%%% N3
\begin{center}
\begin{tikzpicture}[scale=1]
\draw[line width=0.3mm] (-0.2,0) circle (0.2);\draw[line width=0.3mm] (0,0)--(0.5,0);
\draw[line width=0.3mm	] (0.5,0) .. controls (0.8,-0.08) and (1.2,-0.08) .. (1.5,0);
\draw[line width=0.3mm	] (0.5,0) .. controls (1,0.5)  .. (1.5,0);\draw[line width=0.3mm] (1.5,0)--(2,0);
\draw[line width=0.3mm] (2.2,0) circle (0.2);\draw[line width=0.3mm](1,0.4)--(1,1);\draw[line width=0.3mm] (1,1.2) circle (0.2);
	\fill(0,0) circle (0.07);\fill(0.5,0) circle (0.07);\fill(1.5,0) circle (0.07);\fill(2,0) circle (0.07);\fill(1,0.4) circle (0.07);\fill(1,1) circle (0.07);
	{\tiny \node at (01,-0.5){12-3-1-1-1};}
	%%%%%%%%%%%%%%%%%%%%%%
	\draw[line width=0.3mm] (4.8,0) circle (0.2);\draw[line width=0.3mm] (5,0)--(5.5,0);
	\draw[line width=0.3mm	] (5.5,0) .. controls (5.6,-0.09) and (6.1,-0.09) .. (6.2,0);
	\draw[line width=0.3mm	] (5.5,0) .. controls (5.6,0.09) and (6.1,0.09) .. (6.2,0);\draw[line width=0.3mm] (6.2,0)--(6.8,0);
	\draw[line width=0.3mm	] (6.8,0) .. controls (6.9,-0.09) and (7.4,-0.09) .. (7.5,0);
	\draw[line width=0.3mm	] (6.8,0) .. controls (6.9,0.09) and (7.4,0.09) .. (7.5,0);\draw[line width=0.3mm] (7.5,0)--(8,0);
	\draw[line width=0.3mm] (8.2,0) circle (0.2);
	\fill(5,0) circle (0.07);\fill(5.5,0) circle (0.07);\fill(6.2,0) circle (0.07);\fill(6.8,0) circle (0.07);\fill(7.5,0) circle (0.07);\fill(8,0) circle (0.07);
	{\tiny \node at (6.5,-0.5){12-2-2-1-1};}
\end{tikzpicture}
\end{center}
\begin{comment}

%%%% N4
\begin{center}
\begin{tikzpicture}[scale=0.8]
\draw[line width=0.3mm] (-0.2,0) circle (0.2);\draw[line width=0.3mm] (0,0)--(0.5,0);
\draw[line width=0.3mm	] (0.5,0) .. controls (0.6,-0.09) and (1.1,-0.09) .. (1.2,0);
\draw[line width=0.3mm	] (0.5,0) .. controls (0.6,0.09) and (1.1,0.09) .. (1.2,0);\draw[line width=0.3mm] (1.2,0)--(1.8,0);
\draw[line width=0.3mm	] (1.8,0) .. controls (1.9,-0.09) and (2.4,-0.09) .. (2.5,0);
\draw[line width=0.3mm	] (1.8,0) .. controls (1.9,0.09) and (2.4,0.09) .. (2.5,0);\draw[line width=0.3mm] (2.5,0)--(3,0);
\draw[line width=0.3mm] (3.2,0) circle (0.2);
	\fill(0,0) circle (0.07);\fill(0.5,0) circle (0.07);\fill(1.2,0) circle (0.07);\fill(1.8,0) circle (0.07);\fill(2.5,0) circle (0.07);\fill(3,0) circle (0.07);
		{\tiny \node at (1.5,-0.5){12-2-2-1-1};}
\end{tikzpicture}
\end{center}

\end{comment}
%%%% N5
\begin{center}
\begin{tikzpicture}[scale=1]
\draw[line width=0.3mm] (-0.4,-0.3)..controls(-0.6,0) .. (-0.4,0.3);\draw[line width=0.3mm] (-0.4,-0.3)..controls(-0.2,0) .. (-0.4,0.3);
\draw[line width=0.3mm] (-0.4,-0.3)--(0.3,0);\draw[line width=0.3mm] (-0.4,0.3)--(0.3,0);\draw[line width=0.3mm](0.3,0)--(1.3,0);

\draw[line width=0.3mm](1.3,0)--(2,1);\draw[line width=0.3mm](1.3,0)--(2,-1);
\draw[line width=0.3mm] (2.14,1.14) circle (0.2);
\draw[line width=0.3mm] (2.14,-1.14) circle (0.2);
\fill(-0.4,-0.3) circle (0.07);\fill(-0.4,0.3) circle (0.07);\fill(0.3,0) circle (0.07);\fill(1.3,0) circle (0.07);\fill(2,1) circle (0.07);\fill(2,-1) circle (0.07);
{\tiny \node at (0.7,-0.7){11-3-2-1-1};}
%%%% N6
\draw[line width=0.3mm] (3.7,0)--(4.2,0);
\draw[line width=0.3mm] (4.4,0) circle (0.2);
\draw[line width=0.3mm	] (3.7,0) .. controls (4,0.5) and (5,0.5) .. (5.3,0);
\draw[line width=0.3mm	] (3.7,0) .. controls (4,-0.5) and (5,-0.5) .. (5.3,0);
\draw[line width=0.3mm](5.3,0)--(6.3,0);
\draw[line width=0.3mm](6.3,0)--(7,1);\draw[line width=0.3mm](6.3,0)--(7,-1);
\draw[line width=0.3mm] (7.14,1.14) circle (0.2);
\draw[line width=0.3mm] (7.14,-1.14) circle (0.2);
\fill(3.7,0) circle (0.07);\fill(4.2,0) circle (0.07);\fill(3.7,0) circle (0.07);\fill(5.3,0) circle (0.07);\fill(7,1) circle (0.07);\fill(7,-1) circle (0.07);
{\tiny \node at (5.7,-0.7){10-5-1-1-1};}
\end{tikzpicture}
\end{center}
%%%% N6
\begin{comment}

\begin{center}
\begin{tikzpicture}[scale=0.8]
\draw[line width=0.3mm] (-1.3,0)--(-0.8,0);
\draw[line width=0.3mm] (-0.6,0) circle (0.2);
\draw[line width=0.3mm	] (-1.3,0) .. controls (-1,0.5) and (0,0.5) .. (0.3,0);
\draw[line width=0.3mm	] (-1.3,0) .. controls (-1,-0.5) and (0,-0.5) .. (.3,0);
\draw[line width=0.3mm](0.3,0)--(1.3,0);
\draw[line width=0.3mm](1.3,0)--(2,1);\draw[line width=0.3mm](1.3,0)--(2,-1);
\draw[line width=0.3mm] (2.14,1.14) circle (0.2);
\draw[line width=0.3mm] (2.14,-1.14) circle (0.2);
\fill(-1.3,0) circle (0.07);\fill(-0.8,0) circle (0.07);\fill(-1.3,0) circle (0.07);\fill(0.3,0) circle (0.07);\fill(2,1) circle (0.07);\fill(2,-1) circle (0.07);
{\tiny \node at (0.7,-0.7){10-5-1-1-1};}
\end{tikzpicture}
\end{center}

\end{comment}

%%%% N7
\begin{center}
\begin{tikzpicture}[scale=1]
\draw[line width=0.3mm] (-0.4,-0.5)..controls(-0.6,0) .. (-0.4,0.5);\draw[line width=0.3mm] (-0.4,-0.5)..controls(-0.2,0) .. (-0.4,0.5);
\draw[line width=0.3mm] (-0.4,-0.5)--(0.3,-0.5);\draw[line width=0.3mm] (-0.4,0.5)--(0.3,0.5);\draw[line width=0.3mm](0.3,-0.5)--(0.3,0.5);
\draw[line width=0.3mm] (0.3,0.5)--(1.3,0.5);
\draw[line width=0.3mm] (1.5,0.5) circle (0.2);
\draw[line width=0.3mm] (0.3,-0.5)--(1.3,-0.5);
\draw[line width=0.3mm] (1.5,-0.5) circle (0.2);
\fill(-0.4,-0.5) circle (0.07);\fill(-0.4,0.5) circle (0.07);\fill(0.3,0.5) circle (0.07);\fill(1.3,0.5) circle (0.07);\fill(0.3,-0.5) circle (0.07);\fill(1.3,-0.5) circle (0.07);
{\tiny \node at (0.5,-0.9){10-4-2-1-1};}
%%%% N8
\draw[line width=0.3mm] (4.6,-0.3)..controls(4.4,0) .. (4.6,0.3);\draw[line width=0.3mm] (4.6,-0.3)..controls(4.8,0) .. (4.6,0.3);
\draw[line width=0.3mm] (4.6,-0.3)--(5.3,0);\draw[line width=0.3mm] (4.6,0.3)--(5.3,0);\draw[line width=0.3mm](5.3,0)--(6.3,0);
\draw[line width=0.3mm	] (6.3,0) .. controls (6.4,-0.09) and (6.9,-0.09) .. (7,0);
\draw[line width=0.3mm	] (6.3,0) .. controls (6.4, 0.09) and (6.9,0.09) .. (7,0);\draw[line width=0.3mm] (7,0)--(8,0);
\draw[line width=0.3mm] (8.2,0) circle (0.2);
\fill(4.6,-0.3) circle (0.07);\fill (4.6,0.3) circle (0.07);\fill(5.3,0) circle (0.07);\fill(6.3,0) circle (0.07);\fill(7,0) circle (0.07);\fill(8,0) circle (0.07);
{\tiny \node at (6.5,-0.5){10-3-2-2-1};}
\end{tikzpicture}
\end{center}
%%%% N8
\begin{comment}

\begin{center}
\begin{tikzpicture}[scale=0.8]
\draw[line width=0.3mm] (-0.4,-0.3)..controls(-0.6,0) .. (-0.4,0.3);\draw[line width=0.3mm] (-0.4,-0.3)..controls(-0.2,0) .. (-0.4,0.3);
\draw[line width=0.3mm] (-0.4,-0.3)--(0.3,0);\draw[line width=0.3mm] (-0.4,0.3)--(0.3,0);\draw[line width=0.3mm](0.3,0)--(1.3,0);
\draw[line width=0.3mm	] (1.3,0) .. controls (1.4,-0.09) and (1.9,-0.09) .. (2,0);
\draw[line width=0.3mm	] (1.3,0) .. controls (1.4, 0.09) and (1.9,0.09) .. (2,0);\draw[line width=0.3mm] (2,0)--(3,0);
\draw[line width=0.3mm] (3.2,0) circle (0.2);
\fill(-0.4,-0.3) circle (0.07);\fill (-0.4,0.3) circle (0.07);\fill(0.3,0) circle (0.07);\fill(1.3,0) circle (0.07);\fill(2,0) circle (0.07);\fill(3,0) circle (0.07);
{\tiny \node at (1.5,-0.5){10-3-2-2-1};}
\end{tikzpicture}
\end{center}
\end{comment}

%%%% N9
\begin{center}
\begin{tikzpicture}[scale=1]
\draw[line width=0.3mm] (-0.2,0) circle (0.2);\draw[line width=0.3mm] (0,0)--(1,0);
\draw[line width=0.3mm](1,0)--(1.5,0.4);\draw[line width=0.3mm](1,0)--(1.5,-0.4);\draw[line width=0.3mm](1.5,0.4)--(1.5,-0.4);\draw[line width=0.3mm](1.5,-0.4)--(2,0);\draw[line width=0.3mm](1.5,0.4)--(2,0);

\draw[line width=0.3mm] (2,0)--(3,0);
\draw[line width=0.3mm] (3.2,0) circle (0.2);
\fill(0,0) circle (0.07);\fill(1,0) circle (0.07);\fill(1.5,0.4) circle (0.07);\fill(1.5,-0.4) circle (0.07);\fill(2,0) circle (0.07);\fill(3,0) circle (0.07);
{\tiny \node at (1.5,-0.7){10-3-3-1-1};}
%%%% N10
\draw[line width=0.3mm] (4.8,0) circle (0.2);\draw[line width=0.3mm] (5,0)--(6,0);

\draw[line width=0.3mm](6,0)..controls(6.2,0.7) and (7.3,0.7)..(7.5,0);
\draw[line width=0.3mm](6,0)..controls(6.2,-0.7) and (7.3,-0.7)..(7.5,0);
\draw[line width=0.3mm] (7.5,0)--(8.5,0);\draw[line width=0.3mm] (8.7,0) circle (0.2);
\draw[line width=0.3mm] (6.75,-0.5)--(6.75,-0.2);\draw[line width=0.3mm] (6.75,0) circle (0.2);

\fill(5,0) circle (0.07);\fill(6,0) circle (0.07);\fill(7.5,0) circle (0.07);\fill(8.5,0) circle (0.07);\fill(6.75,-0.5) circle (0.07);\fill(6.75,-0.2) circle (0.07);
{\tiny \node at (6.75,-0.9){9-6-1-1-1};}
\end{tikzpicture}
\end{center}
%%%% N10
\begin{comment}

\begin{center}
\begin{tikzpicture}[scale=0.8]
\draw[line width=0.3mm] (-0.2,0) circle (0.2);\draw[line width=0.3mm] (0,0)--(1,0);

\draw[line width=0.3mm](1,0)..controls(1.2,0.7) and (2.3,0.7)..(2.5,0);
\draw[line width=0.3mm](1,0)..controls(1.2,-0.7) and (2.3,-0.7)..(2.5,0);
\draw[line width=0.3mm] (2.5,0)--(3.5,0);\draw[line width=0.3mm] (3.7,0) circle (0.2);
\draw[line width=0.3mm] (1.75,-0.5)--(1.75,-0.2);\draw[line width=0.3mm] (1.75,0) circle (0.2);

\fill(0,0) circle (0.07);\fill(1,0) circle (0.07);\fill(2.5,0) circle (0.07);\fill(3.5,0) circle (0.07);\fill(1.75,-0.5) circle (0.07);\fill(1.75,-0.2) circle (0.07);
{\tiny \node at (1.75,-0.9){9-6-1-1-1};}
\end{tikzpicture}
\end{center}

\end{comment}

%%%% N11
\begin{center}
\begin{tikzpicture}[scale=1]
\draw[line width=0.3mm] (-1.3,0)--(-0.8,0);
\draw[line width=0.3mm] (-0.6,0) circle (0.2);
\draw[line width=0.3mm	] (-1.3,0) .. controls (-1,0.5) and (0,0.5) .. (0.3,0);
\draw[line width=0.3mm	] (-1.3,0) .. controls (-1,-0.5) and (0,-0.5) .. (0.3,0);
\draw[line width=0.3mm](0.3,0)--(1.3,0);
\draw[line width=0.3mm	] (1.3,0) .. controls (1.4,-0.09) and (1.9,-0.09) .. (2,0);
\draw[line width=0.3mm	] (1.3,0) .. controls (1.4, 0.09) and (1.9,0.09) .. (2,0);\draw[line width=0.3mm] (2,0)--(3,0);
\draw[line width=0.3mm] (3.2,0) circle (0.2);
\fill(-1.3,0) circle (0.07);\fill (-0.8,0) circle (0.07);\fill(0.3,0) circle (0.07);\fill(1.3,0) circle (0.07);\fill(2,0) circle (0.07);\fill(3,0) circle (0.07);
{\tiny \node at (1.5,-0.5){9-5-2-1-1};}
%%%% N12
\draw[line width=0.3mm	] (5,-0.4) .. controls (4.8,-0.3) and (4.8,0.3) .. (5,0.4);
\draw[line width=0.3mm	] (5,-0.4) .. controls (5.2,-0.3) and (5.2,0.3) .. (5,0.4);
\draw[line width=0.3mm	] (7,-0.4) .. controls (6.8,-0.3) and (6.8,0.3) .. (7,0.4);
\draw[line width=0.3mm	] (7,-0.4) .. controls (7.2,-0.3) and (7.2,0.3) .. (7,0.4);
\draw[line width=0.3mm	] (5,-0.4)--(7,-0.4); \draw[line width=0.3mm	] (5,0.4)--(7,0.4);

\draw[line width=0.3mm](6,0.4)--(6,0.8);\draw[line width=0.3mm] (6,1) circle (0.2);
\fill(5,-0.4) circle (0.07);\fill (5,0.4) circle (0.07);\fill(7,0.4) circle (0.07);\fill(7,-0.4) circle (0.07);\fill(6,0.4) circle (0.07);\fill(6,0.8) circle (0.07);
{\tiny \node at (6,-0.8){8-5-2-2-1};}
\end{tikzpicture}
\end{center}
%%%% N12
\begin{comment}

\begin{center}
\begin{tikzpicture}[scale=0.8]
\draw[line width=0.3mm	] (0,-0.4) .. controls (-0.2,-0.3) and (-0.2,0.3) .. (0,0.4);
\draw[line width=0.3mm	] (0,-0.4) .. controls (0.2,-0.3) and (0.2,0.3) .. (0,0.4);
\draw[line width=0.3mm	] (2,-0.4) .. controls (1.8,-0.3) and (1.8,0.3) .. (2,0.4);
\draw[line width=0.3mm	] (2,-0.4) .. controls (2.2,-0.3) and (2.2,0.3) .. (2,0.4);
\draw[line width=0.3mm	] (0,-0.4)--(2,-0.4); \draw[line width=0.3mm	] (0,0.4)--(2,0.4);

\draw[line width=0.3mm](1,0.4)--(1,0.8);\draw[line width=0.3mm] (1,1) circle (0.2);
\fill(0,-0.4) circle (0.07);\fill (0,0.4) circle (0.07);\fill(2,0.4) circle (0.07);\fill(2,-0.4) circle (0.07);\fill(1,0.4) circle (0.07);\fill(1,0.8) circle (0.07);
{\tiny \node at (1,-0.8){8-5-2-2-1};}
\end{tikzpicture}
\end{center}

\end{comment}
%%%% N13
\begin{center}
\begin{tikzpicture}[scale=1]

\draw[line width=0.3mm	] (0,-0.4) .. controls (-0.2,-0.3) and (-0.2,0.3) .. (0,0.4);
\draw[line width=0.3mm	] (0,-0.4) .. controls (0.2,-0.3) and (0.2,0.3) .. (0,0.4);
\draw[line width=0.3mm	] (2,-0.4) .. controls (1.8,-0.3) and (1.8,0.3) .. (2,0.4);
\draw[line width=0.3mm	] (2,-0.4) .. controls (2.2,-0.3) and (2.2,0.3) .. (2,0.4);
\draw[line width=0.3mm	] (0,-0.4)--(2,-0.4); \draw[line width=0.3mm	] (0,0.4)--(2,0.4);

\draw[line width=0.3mm](2.13,0)--(2.8,0);\draw[line width=0.3mm] (3,0) circle (0.2);
\fill(0,-0.4) circle (0.07);\fill (0,0.4) circle (0.07);\fill(2,0.4) circle (0.07);\fill(2,-0.4) circle (0.07);\fill(2.13,0) circle (0.07);\fill(2.8,0) circle (0.07);
{\tiny \node at (1,-0.8){8-4-3-2-1};}
%%%% N14
\draw[line width=0.3mm] (4.6,-0.4)..controls(4.4,0) .. (4.6,0.4);\draw[line width=0.3mm] (4.6,-0.4)..controls(4.8,0) .. (4.6,0.4);
\draw[line width=0.3mm] (4.6,-0.4)--(5.3,0);\draw[line width=0.3mm] (4.6,0.4)--(5.3,0);
\draw[line width=0.3mm](5.3,0)--(6.8,0);
\draw[line width=0.3mm] (7.5,0.4)..controls(7.3,0) .. (7.5,-0.4);\draw[line width=0.3mm] (7.5,0.4)..controls(7.7,0) .. (7.5,-0.4);
\draw[line width=0.3mm] (7.5,-0.4)--(6.8,0);\draw[line width=0.3mm] (7.5,0.4)--(6.8,0);
\fill (4.6,-0.4) circle (0.07);\fill (4.6,0.4) circle (0.07);\fill(5.3,0) circle (0.07);\fill(7.5,-0.4) circle (0.07);\fill(6.8,0) circle (0.07);\fill(7.5,0.4) circle (0.07);
{\tiny \node at (6,-0.5){8-3-3-2-2};}
\end{tikzpicture}
\end{center}
%%%% N14
\begin{comment}

\begin{center}
\begin{tikzpicture}[scale=0.8]
\draw[line width=0.3mm] (-0.4,-0.4)..controls(-0.6,0) .. (-0.4,0.4);\draw[line width=0.3mm] (-0.4,-0.4)..controls(-0.2,0) .. (-0.4,0.4);
\draw[line width=0.3mm] (-0.4,-0.4)--(0.3,0);\draw[line width=0.3mm] (-0.4,0.4)--(0.3,0);
\draw[line width=0.3mm](0.3,0)--(1.8,0);
\draw[line width=0.3mm] (2.5,0.4)..controls(2.3,0) .. (2.5,-0.4);\draw[line width=0.3mm] (2.5,0.4)..controls(2.7,0) .. (2.5,-0.4);
\draw[line width=0.3mm] (2.5,-0.4)--(1.8,0);\draw[line width=0.3mm] (2.5,0.4)--(1.8,0);
\fill (-0.4,-0.4) circle (0.07);\fill (-0.4,0.4) circle (0.07);\fill(0.3,0) circle (0.07);\fill(2.5,-0.4) circle (0.07);\fill(1.8,0) circle (0.07);\fill(2.5,0.4) circle (0.07);
{\tiny \node at (1,-0.5){8-3-3-2-2};}
\end{tikzpicture}
\end{center}

\end{comment}
%%%% N15
\begin{center}
\begin{tikzpicture}[scale=1]
\draw[line width=0.3mm	] (0,-0.4) .. controls (-0.2,-0.3) and (-0.2,0.3) .. (0,0.4);
\draw[line width=0.3mm	] (0,-0.4) .. controls (0.2,-0.3) and (0.2,0.3) .. (0,0.4);
\draw[line width=0.3mm	] (2,-0.4)--(2,0.4);
\draw[line width=0.3mm	] (0,-0.4)--(2,-0.4); \draw[line width=0.3mm	] (0,0.4)--(2,0.4);
\draw[line width=0.3mm](2,-0.4)--(2.8,-0.4);\draw[line width=0.3mm] (3,-0.4) circle (0.2);\draw[line width=0.3mm] (2,0.4)--(1.18,0.10);
\draw[line width=0.3mm] (1,0) circle (0.2);

\fill(0,-0.4) circle (0.07);\fill (0,0.4) circle (0.07);\fill(2,0.4) circle (0.07);\fill(2,-0.4) circle (0.07);\fill(1.18,0.10) circle (0.07);\fill(2.8,-0.4) circle (0.07);
{\tiny \node at (1,-0.8){7-7-2-1-1};}

%%%% N16
\draw[line width=0.3mm	] (5,-0.4) .. controls (4.8,-0.3) and (4.8,0.3) .. (5,0.4);
\draw[line width=0.3mm	] (5,-0.4) .. controls (5.2,-0.3) and (5.2,0.3) .. (5,0.4);
\draw[line width=0.3mm	] (7,-0.4)--(7,0.4);
\draw[line width=0.3mm	] (5,-0.4)--(7,-0.4); \draw[line width=0.3mm	] (5,0.4)--(7,0.4);
\draw[line width=0.3mm](7,0.4)--(7.8,0.4);\draw[line width=0.3mm] (8,0.4) circle (0.2);

\draw[line width=0.3mm] (7,-0.4)--(6.18,-0.10);
\draw[line width=0.3mm] (6,0) circle (0.2);

\fill(5,-0.4) circle (0.07);\fill (5,0.4) circle (0.07);\fill(7,0.4) circle (0.07);\fill(7,-0.4) circle (0.07);\fill(6.18,-0.10) circle (0.07);\fill(7.8,0.4) circle (0.07);
{\tiny \node at (6,-0.8){7-7-2-1-1};}

\end{tikzpicture}
\end{center}
%%%% N16
\begin{comment}

\begin{center}
\begin{tikzpicture}[scale=0.8]
\draw[line width=0.3mm	] (0,-0.4) .. controls (-0.2,-0.3) and (-0.2,0.3) .. (0,0.4);
\draw[line width=0.3mm	] (0,-0.4) .. controls (0.2,-0.3) and (0.2,0.3) .. (0,0.4);
\draw[line width=0.3mm	] (2,-0.4)--(2,0.4);
\draw[line width=0.3mm	] (0,-0.4)--(2,-0.4); \draw[line width=0.3mm	] (0,0.4)--(2,0.4);
\draw[line width=0.3mm](2,0.4)--(2.8,0.4);\draw[line width=0.3mm] (3,0.4) circle (0.2);

\draw[line width=0.3mm] (2,-0.4)--(1.18,-0.10);
\draw[line width=0.3mm] (1,0) circle (0.2);

\fill(0,-0.4) circle (0.07);\fill (0,0.4) circle (0.07);\fill(2,0.4) circle (0.07);\fill(2,-0.4) circle (0.07);\fill(1.18,-0.10) circle (0.07);\fill(2.8,0.4) circle (0.07);
{\tiny \node at (1,-0.8){7-7-2-1-1};}
\end{tikzpicture}
\end{center}

\end{comment}
%%%%N17
\begin{center}
\begin{tikzpicture}[scale=1]
\draw[line width=0.3mm] (-1.3,0)--(-0.8,0);
\draw[line width=0.3mm] (-0.6,0) circle (0.2);
\draw[line width=0.3mm	] (-1.3,0) .. controls (-1,0.5) and (0,0.5) .. (0.3,0.5);
\draw[line width=0.3mm	] (-1.3,0) .. controls (-1,-0.5) and (0,-0.5) .. (0.3,-0.5);
\draw[line width=0.3mm	](0.3,0.5)--(0.3,-0.5);\draw[line width=0.3mm	](0.3,0.5)--(1,0);\draw[line width=0.3mm	](1,0)--(0.3,-0.5);
\draw[line width=0.3mm] (1,0)--(1.8,0);
\draw[line width=0.3mm] (2,0) circle (0.2);

\fill(-1.3,0) circle (0.07);\fill (-0.8,0) circle (0.07);\fill(0.3,0.5) circle (0.07);\fill(0.3,-0.5) circle (0.07);\fill(1,0) circle (0.07);\fill(1.8,0) circle (0.07);
{\tiny \node at (0,-0.8){7-6-3-1-1};}
%%%%18

\draw[line width=0.3mm] (3.7,0)--(4.2,0);
\draw[line width=0.3mm] (4.4,0) circle (0.2);
\draw[line width=0.3mm	] (3.7,0) .. controls (4,0.5) and (5,0.5) .. (5.3,0);
\draw[line width=0.3mm	] (3.7,0) .. controls (4,-0.5) and (5,-0.5) .. (5.3,0);
\draw[line width=0.3mm](5.3,0)--(6.3,0);

\draw[line width=0.3mm] (7,0.4)..controls(6.8,0) .. (7,-0.4);\draw[line width=0.3mm] (7,0.4)..controls(7.2,0) .. (7,-0.4);
\draw[line width=0.3mm] (7,-0.4)--(6.3,0);\draw[line width=0.3mm] (7,0.4)--(6.3,0);
\fill (4.3,0) circle (0.07);\fill (4.2,0) circle (0.07);\fill(5.3,0) circle (0.07);\fill(7,-0.4) circle (0.07);\fill(6.3,0) circle (0.07);\fill(7,0.4) circle (0.07);
{\tiny \node at (5.5,-0.7){7-5-3-2-1};}
\end{tikzpicture}
\end{center}

%%%%18
\begin{comment}

\begin{center}
\begin{tikzpicture}[scale=0.8]
\draw[line width=0.3mm] (-1.3,0)--(-0.8,0);
\draw[line width=0.3mm] (-0.6,0) circle (0.2);
\draw[line width=0.3mm	] (-1.3,0) .. controls (-1,0.5) and (0,0.5) .. (0.3,0);
\draw[line width=0.3mm	] (-1.3,0) .. controls (-1,-0.5) and (0,-0.5) .. (.3,0);
\draw[line width=0.3mm](0.3,0)--(1.3,0);

\draw[line width=0.3mm] (2,0.4)..controls(1.8,0) .. (2,-0.4);\draw[line width=0.3mm] (2,0.4)..controls(2.2,0) .. (2,-0.4);
\draw[line width=0.3mm] (2,-0.4)--(1.3,0);\draw[line width=0.3mm] (2,0.4)--(1.3,0);
\fill (-1.3,0) circle (0.07);\fill (-0.8,0) circle (0.07);\fill(0.3,0) circle (0.07);\fill(2,-0.4) circle (0.07);\fill(1.3,0) circle (0.07);\fill(2,0.4) circle (0.07);
{\tiny \node at (0.5,-0.7){7-5-3-2-1};}
\end{tikzpicture}
\end{center}

\end{comment}
%%%%%N19
\begin{center}
\begin{tikzpicture}[scale=1]
\draw[line width=0.3mm](-0.4,0)--(0.4,0);\draw[line width=0.3mm](-0.4,0)--(0,0.4);\draw[line width=0.3mm](-0.4,0)--(0,-0.4);\draw[line width=0.3mm](0,0.4)--(0.4,0);\draw[line width=0.3mm](0,-0.4)--(0.4,0);\draw[line width=0.3mm](0,0.4)--(1.5,0);\draw[line width=0.3mm](0,-0.4)--(1.5,0);
\draw[line width=0.3mm](1.5,0)--(2.4,0);\draw[line width=0.3mm](2.6,0) circle(0.2);
\fill (-0.4,0) circle (0.07);\fill (0.4,0) circle (0.07);\fill(0,-0.4) circle (0.07);\fill(0,0.4) circle (0.07);\fill(1.5,0) circle (0.07);\fill(2.4,0) circle (0.07);
{\tiny \node at (1
	,-0.6){7-4-3-3-1};}

%%%%%N20
\draw[line width=0.3mm] (3.7,0)--(4.2,0);
\draw[line width=0.3mm] (4.4,0) circle (0.2);
\draw[line width=0.3mm	] (3.7,0) .. controls (4,0.5) and (5,0.5) .. (5.3,0.5);
\draw[line width=0.3mm	] (3.7,0) .. controls (4,-0.5) and (5,-0.5) .. (5.3,-0.5);
\draw[line width=0.3mm	](5.3,0.5)--(5.3,-0.5);
\draw[line width=0.3mm	] (6.9,0) .. controls (6.6,0.5) and (5,0.5) .. (5.3,0.5);
\draw[line width=0.3mm	] (6.9,0) .. controls (6.6,-0.5) and (5,-0.5) .. (5.3,-0.5);
\draw[line width=0.3mm] (6.4,0)--(6.9,0);\draw[line width=0.3mm] (6.2,0) circle (0.2);
\fill(3.7,0) circle (0.07);\fill (4.2,0) circle (0.07);\fill(5.3,0.5) circle (0.07);\fill(5.3,-0.5) circle (0.07);\fill(6.9,0) circle (0.07);\fill(6.4,0) circle (0.07);
{\tiny \node at (5.3,-0.8){6-6-4-1-1};}
\end{tikzpicture}
\end{center}
%%%%%N20
\begin{comment}

\begin{center}
\begin{tikzpicture}[scale=0.8]
\draw[line width=0.3mm] (-1.3,0)--(-0.8,0);
\draw[line width=0.3mm] (-0.6,0) circle (0.2);
\draw[line width=0.3mm	] (-1.3,0) .. controls (-1,0.5) and (0,0.5) .. (0.3,0.5);
\draw[line width=0.3mm	] (-1.3,0) .. controls (-1,-0.5) and (0,-0.5) .. (0.3,-0.5);
\draw[line width=0.3mm	](0.3,0.5)--(0.3,-0.5);
\draw[line width=0.3mm	] (1.9,0) .. controls (1.6,0.5) and (0,0.5) .. (0.3,0.5);
\draw[line width=0.3mm	] (1.9,0) .. controls (1.6,-0.5) and (0,-0.5) .. (0.3,-0.5);
\draw[line width=0.3mm] (1.4,0)--(1.9,0);\draw[line width=0.3mm] (1.2,0) circle (0.2);
\fill(-1.3,0) circle (0.07);\fill (-0.8,0) circle (0.07);\fill(0.3,0.5) circle (0.07);\fill(0.3,-0.5) circle (0.07);\fill(1.9,0) circle (0.07);\fill(1.4,0) circle (0.07);
{\tiny \node at (0.3,-0.8){6-6-4-1-1};}
\end{tikzpicture}
\end{center}

\end{comment}
%%%%N21
\begin{center}
\begin{tikzpicture}[scale=1]

\draw[line width=0.3mm	] (0,-0.4) .. controls (-0.2,-0.3) and (-0.2,0.3) .. (0,0.4);
\draw[line width=0.3mm	] (0,-0.4) .. controls (0.2,-0.3) and (0.2,0.3) .. (0,0.4);
\draw[line width=0.3mm	] (2,-0.4) .. controls (1.8,-0.3) and (1.8,0.3) .. (2,0.4);
\draw[line width=0.3mm	] (2,-0.4) .. controls (2.2,-0.3) and (2.2,0.3) .. (2,0.4);
\draw[line width=0.3mm	] (0,-0.4)--(2,-0.4); \draw[line width=0.3mm	] (0,0.4)--(0.6,0.4);\draw[line width=0.3mm	] (1.4,0.4)--(2,0.4);
\draw[line width=0.3mm	] (0.6,0.4) .. controls (0.8,0.6) and (1.2,0.6) .. (1.4,0.4);
\draw[line width=0.3mm	] (0.6,0.4) .. controls (0.8,0.2) and (1.2,0.2) .. (1.4,0.4);
\fill(0,-0.4) circle (0.07);\fill (0,0.4) circle (0.07);\fill (2,-0.4) circle (0.07);\fill(2,0.4) circle (0.07);\fill(0.6,0.4) circle (0.07);\fill(1.4,0.4) circle (0.07);
{\tiny \node at (1,-0.7){6-6-2-2-2};}

%%%%N22
\draw[line width=0.3mm] (3.7,0)--(4.2,0);
\draw[line width=0.3mm] (4.4,0) circle (0.2);
\draw[line width=0.3mm	] (3.7,0) .. controls (4,0.5) and (5,0.5) .. (5.3,0);
\draw[line width=0.3mm	] (3.7,0) .. controls (4,-0.5) and (5,-0.5) .. (5.3,0);
\draw[line width=0.3mm](5.3,0)--(6.3,0);

\draw[line width=0.3mm	] (6.3,0) .. controls (6.6,0.5) and (7.6,0.5) .. (7.9,0);
\draw[line width=0.3mm	] (6.3,0) .. controls (6.6,-0.5) and (7.6,-0.5) .. (7.9,0);

\draw[line width=0.3mm] (7.9,0)--(7.4,0);
\draw[line width=0.3mm] (7.2,0) circle (0.2);

\fill (3.7,0) circle (0.07);\fill(4.2,0) circle (0.07);\fill (7.4,0) circle (0.07);\fill(5.3,0) circle (0.07);\fill(6.3,0) circle (0.07);\fill(7.9,0) circle (0.07);

{\tiny \node at (5.8,-0.4){6-5-5-1-1};}
\end{tikzpicture}
\end{center}

%%%%N22
\begin{comment}

\begin{center}
\begin{tikzpicture}[scale=0.8]
\draw[line width=0.3mm] (-1.3,0)--(-0.8,0);
\draw[line width=0.3mm] (-0.6,0) circle (0.2);
\draw[line width=0.3mm	] (-1.3,0) .. controls (-1,0.5) and (0,0.5) .. (0.3,0);
\draw[line width=0.3mm	] (-1.3,0) .. controls (-1,-0.5) and (0,-0.5) .. (0.3,0);
\draw[line width=0.3mm](0.3,0)--(1.3,0);

\draw[line width=0.3mm	] (1.3,0) .. controls (1.6,0.5) and (2.6,0.5) .. (2.9,0);
\draw[line width=0.3mm	] (1.3,0) .. controls (1.6,-0.5) and (2.6,-0.5) .. (2.9,0);

\draw[line width=0.3mm] (2.9,0)--(2.4,0);
\draw[line width=0.3mm] (2.2,0) circle (0.2);

\fill (-1.3,0) circle (0.07);\fill(-0.8,0) circle (0.07);\fill (2.4,0) circle (0.07);\fill(0.3,0) circle (0.07);\fill(1.3,0) circle (0.07);\fill(2.9,0) circle (0.07);

{\tiny \node at (0.8,-0.4){6-5-5-1-1};}
\end{tikzpicture}
\end{center}

\end{comment}
%%%%N23
\begin{center}
\begin{tikzpicture}[scale=1]

\draw[line width=0.3mm	] (0,-0.4) .. controls (-0.2,-0.3) and (-0.2,0.3) .. (0,0.4);
\draw[line width=0.3mm	] (0,-0.4) .. controls (0.2,-0.3) and (0.2,0.3) .. (0,0.4);
\draw[line width=0.3mm	] (0.5,-0.4)--(0.5,0.4);
\draw[line width=0.3mm	] (0,-0.4)--(0.5,-0.4); \draw[line width=0.3mm] (0,0.4)--(0.5,0.4);
\draw[line width=0.3mm] (2,0)--(0.5,0.4);\draw[line width=0.3mm] (2,0)--(0.5,-0.4);
\draw[line width=0.3mm] (2,0)--(1.2,0);
\draw[line width=0.3mm] (1.02,0) circle (0.15);
\fill (0,-0.4) circle (0.07);\fill(0,0.4) circle (0.07);\fill(0.5,-0.4) circle (0.07);\fill(0.5,0.4) circle (0.07);\fill(2,0) circle (0.07);\fill(1.2,0) circle (0.07);

{\tiny \node at (0.8,-0.7){6-5-4-2-1};}
%%%%N24
\draw[line width=0.3mm	] (5,-0.4) .. controls (4.8,-0.3) and (4.8,0.3) .. (5,0.4);
\draw[line width=0.3mm	] (5,-0.4) .. controls (5.2,-0.3) and (5.2,0.3) .. (5,0.4);
\draw[line width=0.3mm	] (7,-0.4) .. controls (6.8,-0.3) and (6.8,0.3) .. (7,0.4);
\draw[line width=0.3mm	] (7,-0.4) .. controls (7.2,-0.3) and (7.2,0.3) .. (7,0.4);
\draw[line width=0.3mm	] (5,-0.4)--(7,-0.4); \draw[line width=0.3mm	] (5,0.4)--(7,0.4);\draw[line width=0.3mm	] (6,0.4)--(6,-0.4);

\fill(5,-0.4) circle (0.07);\fill (5,0.4) circle (0.07);\fill (7,-0.4) circle (0.07);\fill(7,0.4) circle (0.07);\fill(6,-0.4) circle (0.07);\fill(6,0.4) circle (0.07);
{\tiny \node at (6,-0.7){6-4-4-2-2};}
\end{tikzpicture}
\end{center}
%%%%N24
\begin{comment}

\begin{center}
\begin{tikzpicture}[scale=0.8]
\draw[line width=0.3mm	] (0,-0.4) .. controls (-0.2,-0.3) and (-0.2,0.3) .. (0,0.4);
\draw[line width=0.3mm	] (0,-0.4) .. controls (0.2,-0.3) and (0.2,0.3) .. (0,0.4);
\draw[line width=0.3mm	] (2,-0.4) .. controls (1.8,-0.3) and (1.8,0.3) .. (2,0.4);
\draw[line width=0.3mm	] (2,-0.4) .. controls (2.2,-0.3) and (2.2,0.3) .. (2,0.4);
\draw[line width=0.3mm	] (0,-0.4)--(2,-0.4); \draw[line width=0.3mm	] (0,0.4)--(2,0.4);\draw[line width=0.3mm	] (1,0.4)--(1,-0.4);

\fill(0,-0.4) circle (0.07);\fill (0,0.4) circle (0.07);\fill (2,-0.4) circle (0.07);\fill(2,0.4) circle (0.07);\fill(1,-0.4) circle (0.07);\fill(1,0.4) circle (0.07);
{\tiny \node at (1,-0.7){6-4-4-2-2};}
\end{tikzpicture}
\end{center}

\end{comment}
%%%%N25
\begin{center}
\begin{tikzpicture}[scale=1]

\draw[line width=0.3mm	] (0,-0.4) .. controls (-0.2,-0.3) and (-0.2,0.3) .. (0,0.4);
\draw[line width=0.3mm	] (0,-0.4) .. controls (0.2,-0.3) and (0.2,0.3) .. (0,0.4);

\draw[line width=0.3mm	] (0,-0.4)--(1.4,-0.4); \draw[line width=0.3mm	] (0,0.4)--(1.4,0.4);\draw[line width=0.3mm	] (1.4,0.4)--(1.1,0);
\draw[line width=0.3mm	] (1.4,0.4)--(1.7,0);\draw[line width=0.3mm	] (1.7,0)--(1.1,0);\draw[line width=0.3mm	] (1.4,-0.4)--(1.1,0);\draw[line width=0.3mm	] (1.4,-0.4)--(1.7,0);
\fill(0,-0.4) circle (0.07);\fill (0,0.4) circle (0.07);\fill (1.4,-0.4) circle (0.07);\fill(1.4,0.4) circle (0.07);\fill(1.1,0) circle (0.07);\fill(1.7,0) circle (0.07);
{\tiny \node at (0.7,-0.7){5-5-3-3-2};}
%%%%N26
	\draw[line width=0.3mm	] (5,-0.4)--(5,0.4);\draw[line width=0.3mm](5,0.4)--(5.5,0);\draw[line width=0.3mm](5,-0.4)--(5.5,0);
\draw[line width=0.3mm	] (5.5,0)--(6.5,0);	\draw[line width=0.3mm](6.5,0)--(7,0.4);   \draw[line width=0.3mm](6.5,0)--(7,-0.4);
\draw[line width=0.3mm	] (7,0.4)--(7,-0.4);\draw[line width=0.3mm](5,0.4)--(7,0.4);\draw[line width=0.3mm](5,-0.4)--(7,-0.4);
\fill(5,-0.4) circle (0.07);\fill (5,0.4) circle (0.07);\fill (5.5,0) circle (0.07);\fill(6.5,0) circle (0.07);\fill(7,-0.4) circle (0.07);\fill(7,0.4)circle (0.07);
{\tiny \node at (6,-0.7){4-4-4-3-3};}
\end{tikzpicture}
\end{center}
%%%%N26
\begin{comment}

\begin{center}
	\begin{tikzpicture}[scale=0.8]
	\draw[line width=0.3mm	] (0,-0.4)--(0,0.4);\draw[line width=0.3mm](0,0.4)--(0.5,0);\draw[line width=0.3mm](0,-0.4)--(0.5,0);
	\draw[line width=0.3mm	] (0.5,0)--(1.5,0);	\draw[line width=0.3mm](1.5,0)--(2,0.4);   \draw[line width=0.3mm](1.5,0)--(2,-0.4);
	\draw[line width=0.3mm	] (2,0.4)--(2,-0.4);\draw[line width=0.3mm](0,0.4)--(2,0.4);\draw[line width=0.3mm](0,-0.4)--(2,-0.4);
	\fill(0,-0.4) circle (0.07);\fill (0,0.4) circle (0.07);\fill (0.5,0) circle (0.07);\fill(1.5,0) circle (0.07);\fill(2,-0.4) circle (0.07);\fill(2,0.4)circle (0.07);
	{\tiny \node at (1,-0.7){4-4-4-3-3};}
	\end{tikzpicture}
	
\end{center}

\end{comment}

The functional invariants of the above K3 elliptic surfaces are given the table below. We notice that they are all defined over $\mathbb Q$ except for the two cases corresponding to the cusp split 7-7-2-1-1 which are defined over the quadratic imaginary field ${\mathbb Q}(\sqrt{-7})$. These invariant are computed as genus zero Beliy maps of the Riemann surfaces $X(\G)$ using the maple program  in \cite{compute}. These Beliy maps can be looked at as the J-invariants of the elliptic surfaces, that is the J-invariant of the elliptic curves over the function field of the base curve, or as a way to express the elliptic modular function $j$ as a rational function of a chosen Hauptmodul of the genus zero modular subgroup.

{\tiny
	%\fontsize{6}{8}\selectfont
	\begin{tabular}{|l|c|}
		\hline
		Cusp split& $J-$function\\
		\hline
14-1-1-1-1& $\displaystyle%
%{\frac {1}{442368}}\,
{\frac { \left( 2197\,{x}^{6}+384+2704\,{x}^{4}+2080\,{x}^{2} \right) ^{3}}{169\,{x}^{4}+169\,{x}^{2}+128}}$\\
	\hline
13-2-1-1-1 & $\displaystyle%
	%{\frac {1}{2616440389632}}\,
	{\frac { \left( 47045881\,{x}^{6}+69330772\,{x}^{5}
			+28670620\,{x}^{4}-1536416\,{x}^{3}-2223760\,{x}^{2}-10944\,x+10816 \right) ^{3}}{{x}^{2} \left( 6859\,{x}^{3}+6859\,{x}^{2}+912\,x-676 \right) 	}}$\\
		\hline
	12-3-1-1-1&	 $\displaystyle%
		%{\frac {1}{108}}\,
		{\frac { \left( 16\,{x}^{6}+1+16\,{x}^{3} \right) ^{3}}{{x}^{3} \left( x+1 \right)  \left( {x}^{2}-x+1 \right) }}$\\
		\hline
	12-2-2-1-1&	 $\displaystyle%
		%{\frac {4}{27}}\,
		{\frac { \left( {x}^{6}+3+6\,{x}^{4}+9\,{x}^{2} \right) ^{3}}{ \left( {x}^{2}+1 \right) ^{2} \left( {x}^{2}+4 \right) 	}}$\\
		\hline
		11-3-2-1-1&	 $\displaystyle%
			%{\frac {1}{7247757312}}\,
			{\frac { \left( {x}^{6}+54\,{x}^{5}+1023\,{x}^{4}+6644\,{x}^{3}+495\,{x}^{2}
					-15114\,x+14641 \right) ^{3}}{{x}^{3} \left( x-1 \right) ^{2} \left( 3\,{x}^{2}+118\,x+1331 \right) }}$\\
				\hline
			10-5-1-1-1&	 $\displaystyle%
				% -{\frac {1}{927712935936}}\,
				{\frac { \left( 16777216\,{x}^{6}+23068672\,{x}^{5}
						+11468800\,{x}^{4}+2048000\,{x}^{3}+44800\,{x}^{2}-672\,x+1 \right) ^{3}}{{x}^{5} \left( 2048\,{x}^{3}+2048\,{x}^{2}+616\,x-1 \right) 
						}}$\\
					\hline
				10-4-2-1-1&		 $\displaystyle%
					% {\frac {1}{27648}}\,
				{\frac { \left( {x}^{6}+4\,{x}^{5}+20\,{x}^{4}+64-64\,x \right) ^{3}}{{x}^{4} \left( x-1 \right) ^{2}
							 \left( {x}^{2}+4\,x+20 \right) }}$\\
						 \hline
			10-3-2-2-1&			 $\displaystyle%
						%-{\frac {1}{442368}}\,
				{\frac { \left( {x}^{6}-384+48\,{x}^{4}+288\,{x}^{2} \right) ^{3}}{ \left( {x}^{2}+1 \right) ^{3}
								 \left( 3\,{x}^{2}+128 \right) }}$\\ \hline
				10-3-3-1-1&				 $\displaystyle%
							%{\frac {1}{1728}}\,
						{\frac { \left( 16384\,{x}^{6}-36864\,{x}^{5}+26880\,{x}^{4}-6400\,{x}^{3}
									+1+24\,x \right) ^{3}}{{x}^{3} \left( x-1 \right)  \left( 24\,x+1 \right) ^{2} \left( 8\,x-5 \right) ^{2}}}$\\
								\hline
							9-5-2-1-1&	 $\displaystyle%
								%-{\frac {1}{12230590464}}\,
							{\frac { \left( 6561\,{x}^{6}+8748\,{x}^{5}+20412\,{x}^{4}
										+6048\,{x}^{3}+3312\,{x}^{2}-1344\,x+64 \right) ^{3}}{{x}^{6} \left( 27\,{x}^{3}+27\,{x}^{2}+72\,x-4 \right) }}$\\
									\hline
						9-6-1-1-1&				 $\displaystyle%
									%{\frac {1}{7247757312}}\,
								{\frac { \left( 59049\,{x}^{6}-223074\,{x}^{5}+273375\,{x}^{4}
											-114300\,{x}^{3}+6295\,{x}^{2}+254\,x+1 \right) ^{3}}{{x}^{5} \left( x-1 \right) ^{2} \left( 135\,{x}^{2}-234\,x-1 \right) }}$\\
										\hline
										
							8-5-2-2-1&				 $\displaystyle%
											%{\frac {1}{6912}}\,
									{\frac { \left( 1024\,{x}^{6}-1024\,{x}^{5}+1+20\,{x}^{2}+4\,x \right) ^{3}}{{x}^{5} \left( x-1 \right) 
													 \left( 20\,{x}^{2}+4\,x+1 \right) ^{2}}}$\\
												 \hline
						8-4-3-2-1&			 $\displaystyle%
												% {\frac {1}{229582512}}\,
												{\frac { \left( 256\,{x}^{6}-1536\,{x}^{5}+2112\,{x}^{4}+64\,{x}^{3}
														-615\,{x}^{2}+192\,x+256 \right) ^{3}}{{x}^{4} \left( x-1 \right) ^{3} \left( 2\,x+1 \right) ^{2} \left( x-4 \right) }}$\\
													\hline
								8-3-3-2-2&					 $\displaystyle%
													%{\frac {1}{108}}\,
													{\frac { \left( 16\,{x}^{6}+12+48\,{x}^{4}+45\,{x}^{2} \right) ^{3}}{ \left( {x}^{2}+1 \right) ^{3}
															 \left( 3\,{x}^{2}+4 \right) ^{2}}}$\\
														 \hline

				%98\,
			7-7-2-1-1&		$\displaystyle%
			\frac{ 
					% \left( -{\frac {203}{294912}}\,\xi-{\frac {4459}{884736}}	 \right) 1/14\\
					\left(14 {x}^{6}+ 7\left( -45+7\,\xi \right) {x}^{4}+ 7\left( -29-9\,\xi \right) {x}^{2}+39-21\,\xi \right) ^{3}}{ \left( {x}^{2}+1 \right) ^{7} \left( 98\,{x}^{2}+49+13\,\xi \right) },\	{\xi}^{2}+7=0$\\
				\hline
			
					%{\frac {1}{452984832}}\,
				7-6-3-1-1&	$\displaystyle	{\frac { \left( 729\,{x}^{6}-5508\,{x}^{5}+10620\,{x}^{4}-5088\,{x}^{3}
							+624\,{x}^{2}+960\,x+64 \right) ^{3}}{{x}^{6} \left( x-1 \right) ^{3} \left( 81\,{x}^{2}-396\,x-28 \right) }}$\\
						\hline
				
						%{\frac {1}{2295825120000000000}}\,
					7-5-3-2-1&	$\displaystyle	{\frac { \left( 1073741824\,{x}^{6}
								-3271557120\,{x}^{5}+3379200000\,{x}^{4}-1254800000\,{x}^{3}+71700000\,{x}^{2}+3450000\,x
								+15625 \right) ^{3}}{{x}^{5} \left( x-1 \right) ^{3} \left( 32\,x-35 \right) ^{2} \left( 1024\,x+5 \right) }}$\\
							\hline
						
							%{\frac {1}{78732}}\,
					7-4-3-3-1&		$\displaystyle	{\frac { \left( 65536\,{x}^{6}-49152\,{x}^{5}-16896\,{x}^{4}+1472\,{x}^{3}
									+1185\,{x}^{2}+240\,x+16 \right) ^{3}}{{x}^{4} \left( x-1 \right)  \left( 256\,{x}^{2}+80\,x+7 \right) ^{3}}}$\\
								\hline
						
								% -{\frac {4}{27}}\,
						6-6-4-1-1&		$\displaystyle	{\frac { \left( {x}^{6}-3+3\,{x}^{4} \right) ^{3}}{ \left( {x}^{2}+1 \right) ^{6} \left( 3\,{x}^{2}+4 \right) }}$\\
						\hline
						
								%{\frac {4}{27}}\,
								6-6-2-2-2&		$\displaystyle{\frac { \left( {x}^{6}+1+{x}^{3} \right) ^{3}}{{x}^{6} \left( x+1 \right) ^{2} \left( {x}^{2}-x+1 \right) ^{2}}}$\\
								\hline
						
								%{\frac {1}{442368}}\,
						6-5-5-1-1&	$\displaystyle	{\frac { \left( 625\,{x}^{6}+384+1680\,{x}^{4}+1440\,{x}^{2} \right) ^{3}}{ \left( {x}^{2}+1 \right) ^{5}
										 \left( 125\,{x}^{2}+128 \right) }}$\\
									 \hline
						
								6-5-4-2-1&		$\displaystyle
								% {\frac {1}{16875000000}}\,
									{\frac { \left( 256\,{x}^{6}+3840\,{x}^{5}+15625-25000\,{x}^{3}
											+65625\,{x}^{2}-56250\,x \right) ^{3}}{{x}^{5} \left( x-1 \right) ^{4} \left( 8\,x-5 \right) ^{2} \left( x+15 \right) }}$\\
										\hline
									
										% -{\frac {1}{108}}\,
									6-4-4-2-2&		$\displaystyle {\frac { \left( {x}^{6}-12+6\,{x}^{4}+9\,{x}^{2} \right) ^{3}}{ \left( {x}^{2}+1 \right) ^{4} \left( {x}^{2}+4 \right) ^{2}
												}}$\\
											\hline
									
											% -{\frac {1}{19683}}\,
										5-5-3-3-2&	$\displaystyle	{\frac{ \left( 625\,{x}^{6}+81-13005\,{x}^{4}+2835\,{x}^{2} \right) ^{3}}{ \left( {x}^{2}+1 \right) ^{5}
													 \left( 125\,{x}^{2}-3 \right) ^{3}}}$\\
												 \hline
									
												% -{\frac {1}{108}}\,
										4-4-4-3-3&		$\displaystyle	{\frac { \left( {x}^{6}+1-14\,{x}^{3} \right) ^{3}}{{x}^{3} \left( 1+x \right) ^{4} \left( {x}^{2}-x+1 \right) ^{4}
														}}$\\
													\hline
										\end{tabular}

}%%%tiny

\section{Further developments: The genus 1 case}
It would be  interesting to study the case of the genus 1 torsion-free subgroups of $\PS$. Their modular curves are thus smooth curves of genus 1. The Riemann-Hurwitz formula takes the form
\begin{equation}\label{h-r-g1}
\mu=6h.
\end{equation}
Similarly to the genus zero case, these groups correspond to rooted trivalent maps on the torus while the conjugacy classes correspond to unrooted trivalent maps on the torus. This correspondence is again done through the permutation groups.
 It should be understood that a trivalent map on the torus is a map that is not homotopic to a point, in other words,  it does not correspond to a map on a sphere. In particular, such a map should have $2h$ vertices,  $3h$ edges and $h$ faces. 
 
 Maps on the torus have been studied and enumerated in many cases by several authors including \cite{kr-om1,kr-om2, me-ne,wa2,wa-le}. An explicit closed formula for the number of genus 1 torsion-free modular subgroups is not yet available and is the topic of  future work, however, one can compute numerically their number for any given index by working with
 the permutation groups. Our calculation turns out to yield the same information as the one in Table 3 and Table 4 of \cite{kr-om1}.

If we denote by $N(\mu)$ the number of genus 1 torsion-free subgroups and $N^+(\mu)$ the number of their conjugacy classes, we have the following table:

\begin{center}
	\begin{tabular}{l|l|l|}
		$\mu$&$N(\mu)$&$N^+(\mu)$\\
		\hline
		6&1&1\\
		12&28&5\\
		18&664&46\\
		24&14912&669\\
		30&326496&11096\\
		36&7048192&19688\\
		\hline
	\end{tabular}
	\end{center}

Below are the maps for the conjugacy classes of genus 1 torsion-free modular subgroups at index 6 and 12 together with their cusp splits. In particular, the map at index 6 has only one face with 6 edges surrounding it (each edge is counting twice, or by considering all edges as double edges with opposite directions). 
	
	
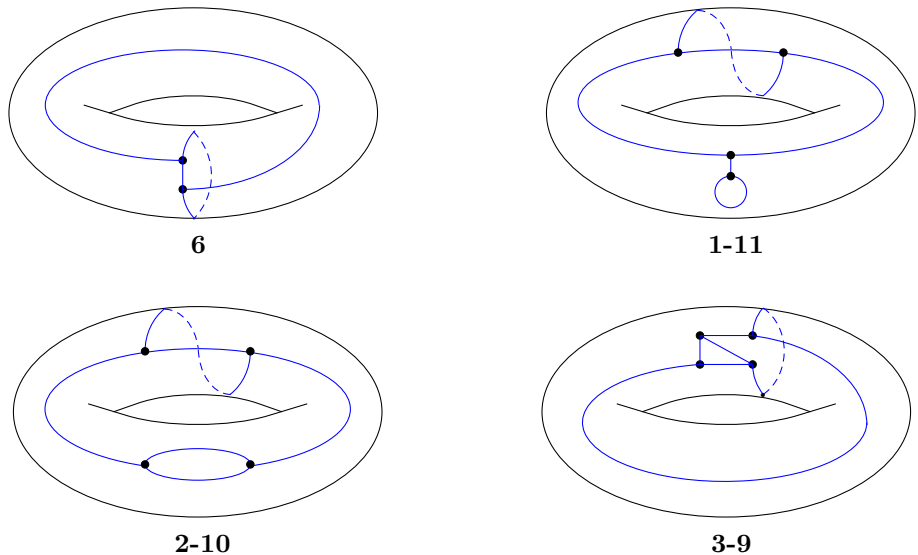
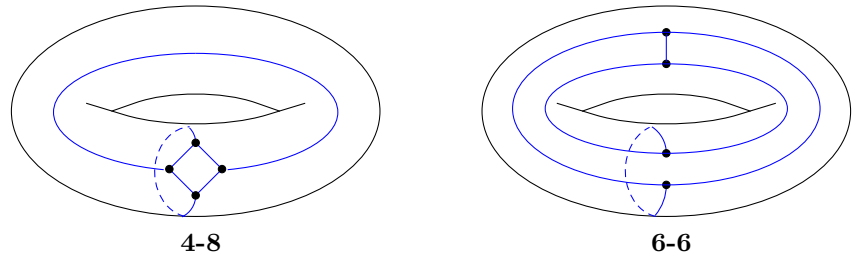
\begin{center}
		\begin{tikzpicture}[scale=0.7]

		%Hole
		\begin{scope}[scale=0.8]
		\path[rounded corners=24pt] (-2,0)--(0,.9)--(2,0) (-1.3,0)--(0,-1)--(1.3,0);
		\draw[rounded corners=28pt] (-2.6,.19)--(0,-.6)--(2.6,.19);
		\draw[rounded corners=24pt] (-2,0)--(0,.7)--(2,0);
		\end{scope}
		
		%Cut 2
		\draw[] (0,0) ellipse (3.5 and 2);

		\draw[fill=black] (-.2,-1.45) circle (2pt);
		\draw[fill=black] (-.2,-.9) circle (2pt);
		
		\draw[blue] (-.2,-.9) arc (180: 135 : 0.8 );
		\draw[blue] (-.2,-.9) -- (-.2,-1.45);
		\draw[blue] (-.2,-1.45) arc (180: 225 : 0.8 );
		\draw[densely dashed,blue] (0,-.35) arc (45: -45 : 1.18 );
		\draw[blue] (-.2,-1.45) arc (-90: 0 : 2.6 and 1.6 );
		\draw[blue] (2.39,0.15) arc (0: 270 : 2.6 and 1.05 );
		
		\node at (0,-2.5) {$\textbf{
				6}$};
		
		\end{tikzpicture}
		\hspace{2cm}	%\end{center}
	%	\begin{center}
	\begin{tikzpicture}[scale=0.7]

%Hole
\begin{scope}[scale=0.8]
\path[rounded corners=24pt] (-2,0)--(0,.9)--(2,0) (-1.3,0)--(0,-1)--(1.3,0);
\draw[rounded corners=28pt] (-2.6,.19)--(0,-.6)--(2.6,.19);
\draw[rounded corners=24pt] (-2,0)--(0,.7)--(2,0);
\end{scope}
%Cut 1
\draw[blue] (0,-.8) arc (270:90:2.9 and 1);
\draw[blue] (0,-.8) arc (-90:90:2.9 and 1);

%Cut 2
\draw[] (0,0) ellipse (3.5 and 2);
\node at (0,-.8) {$\textbf{.}$};
\draw[blue] (0,-.8) -- (0,-1.2);
\node at (0,-1.2) {$\textbf{.}$};
\draw[blue] (0,-1.5) ellipse (0.3 and 0.3);
\node at (-1,1.15) {$\textbf{.}$};
\node at (1,1.15) {$\textbf{.}$};
\draw[blue] (-1,1.15) arc (180:125:0.9 and 1);
\draw[blue] (1,1.15) arc (180:125:-0.9 and -1);
%\node at (0.63,0.33) {$\textbf{.}$};
\draw[densely dashed,blue] (0.63,0.33) arc (-90:0:-0.62 and 0.9);
%\node at (-0.63,1.95) {$\textbf{.}$};
\draw[densely dashed,blue] (-0.63,1.95) arc (-90:-10:0.65 and -0.9);

\draw[fill=black] (0,-.8) circle (2pt);
\draw[fill=black] (0,-1.2) circle (2pt);
\draw[fill=black] (1,1.15) circle (2pt);
\draw[fill=black] (-1,1.15) circle (2pt);

\node at (0,-2.5) {$\textbf{
		1-11}$};

\end{tikzpicture}
\end{center}
\begin{center}
		\begin{tikzpicture}[scale=0.7]

	%Hole
	\begin{scope}[scale=0.8]
	\path[rounded corners=24pt] (-2,0)--(0,.9)--(2,0) (-1.3,0)--(0,-1)--(1.3,0);
	\draw[rounded corners=28pt] (-2.6,.19)--(0,-.6)--(2.6,.19);
	\draw[rounded corners=24pt] (-2,0)--(0,.7)--(2,0);
	\end{scope}
	%Cut 1
	\draw[blue] (0,1.2) arc (270:110:2.9 and -1.15);
	\draw[blue] (0,1.2) arc (-90:70:2.9 and -1.15);
	
	%Cut 2
	\draw[] (0,0) ellipse (3.5 and 2);
	\node at (-1,1.15) {$\textbf{.}$};
	\node at (1,1.15) {$\textbf{.}$};
	\node at (-1,-1) {$\textbf{.}$};
	\node at (1,-1) {$\textbf{.}$};
	
	\draw[blue] (-1,-1) arc (180:0:1 and 0.3);
	\draw[blue] (-1,-1) arc (180:0:1 and -0.3);
	
	\draw[blue] (-1,1.15) arc (180:125:0.9 and 1);
	\draw[blue] (1,1.15) arc (180:125:-0.9 and -1);
	%\node at (0.63,0.33) {$\textbf{.}$};
	\draw[densely dashed,blue] (0.63,0.33) arc (-90:0:-0.62 and 0.9);
	%\node at (-0.63,1.95) {$\textbf{.}$};
	\draw[densely dashed,blue] (-0.63,1.95) arc (-90:-10:0.65 and -0.9);

	\draw[fill=black] (1,1.15) circle (2pt);
	\draw[fill=black] (-1,1.15) circle (2pt);
	\draw[fill=black] (-1,-1) circle (2pt);
	\draw[fill=black] (1,-1) circle (2pt);
	
	\node at (0,-2.5) {$\textbf{
			2-10}$};
		\end{tikzpicture}\hspace{2cm}	
%\end{center}
%\begin{center}
	\begin{tikzpicture}[scale=0.7]

	%Hole
	\begin{scope}[scale=0.8]
	\path[rounded corners=24pt] (-2,0)--(0,.9)--(2,0) (-1.3,0)--(0,-1)--(1.3,0);
	\draw[rounded corners=28pt] (-2.6,.19)--(0,-.6)--(2.6,.19);
	\draw[rounded corners=24pt] (-2,0)--(0,.7)--(2,0);
	\end{scope}
	
	%Cut 2
	\draw[] (0,0) ellipse (3.5 and 2);

	\draw[fill=black] (-.5,1.45) circle (2pt);
	\draw[fill=black] (-.5,.9) circle (2pt);
	
	\draw[fill=black] (.5,1.45) circle (2pt);
	\draw[fill=black] (.5,.9) circle (2pt);
	
	\draw[blue] (-.5,1.45) -- (-.5,.9);
	\draw[blue] (.5,1.45) -- (-.5,1.45);
	\draw[blue] (-.5,.9) -- (.5,.9);
	\draw[blue] (-.5,1.45) -- (.5,.9);
	
	\draw[blue] (-.5,.9) arc (80: -180 :-2.7 and 1.12);
	
	\draw[blue] (2.67,-0.25) arc (0: 80 :2.7 and 1.73);

	\draw[blue] (.5,.9) arc (180: 220 : 0.9  );
	
	\draw[blue] (.5,1.45) arc (180: 140 : 0.8  );
	
	\node at (0.7,.33) {$\textbf{.}$};
	
	\draw[densely dashed,blue] (0.7,.33) arc (-60: 60 :.8 and 0.95);
	
	\node at (0,-2.5) {$\textbf{
			3-9}$};

	\end{tikzpicture}
\end{center}
\begin{center}
\begin{tikzpicture} [scale=0.7]

%Hole
\begin{scope}[scale=0.8]
\path[rounded corners=24pt] (-2,0)--(0,.9)--(2,0) (-1.3,0)--(0,-1)--(1.3,0);
\draw[rounded corners=28pt] (-2.6,.19)--(0,-.6)--(2.6,.19);
\draw[rounded corners=24pt] (-2,0)--(0,.7)--(2,0);
\end{scope}

%Cut 2
\draw[] (0,0) ellipse (3.5 and 2);
\draw[blue] (0,1.1) arc (90:-77:-2.7 and 1.12);
\draw[blue] (0,1.1) arc (90:-77:2.7 and 1.12);
\draw[blue] (0,-1.6) arc (0:-65: 0.45);
\draw[densely dashed,blue] (-0.2,-2) arc (260:93: 0.7 and 0.86);
\draw[blue] (0,-.6) arc (0:45: 0.4);

\draw[blue] (0,-0.6) -- (-0.5,-1.1) -- (0,-1.6) -- (0.5,-1.1) -- (0,-0.6);

\draw[fill=black] (0,-0.6) circle (2pt);
\draw[fill=black] (-0.5,-1.1) circle (2pt);
\draw[fill=black] (0.5,-1.1) circle (2pt);
\draw[fill=black] (0,-1.6) circle (2pt);

\node at (0,-0.6) {$\textbf{.}$};
\node at (-0.5,-1.1) {$\textbf{.}$};
\node at (0.5,-1.1) {$\textbf{.}$};
\node at (0,-1.6) {$\textbf{.}$};

\node at (0,-2.5) {$\textbf{
		4-8}$};

\end{tikzpicture}
%\end{center}
%\begin{center}
		\begin{tikzpicture}[scale=0.7]

	%Hole
	\begin{scope}[scale=0.8]
	\path[rounded corners=24pt] (-2,0)--(0,.9)--(2,0) (-1.3,0)--(0,-1)--(1.3,0);
	\draw[rounded corners=28pt] (-2.6,.19)--(0,-.6)--(2.6,.19);
	\draw[rounded corners=24pt] (-2,0)--(0,.7)--(2,0);
	\end{scope}
	
	%Cut 2
	\draw[] (0,0) ellipse (3.5 and 2);
	\draw[blue] (0,-1.4) arc (0:-40: .9);
	\draw[densely dashed,blue] (-0.2,-2) arc (260:110 : 0.7 and 0.88);
	\draw[blue] (0,-.8) arc (0: 57 : 0.6);

	\draw[fill=black] (0,-1.4) circle (2pt);
	\draw[fill=black] (0,-.8) circle (2pt);
	
	\draw[fill=black] (0,1.5) circle (2pt);
	\draw[fill=black] (0,.9) circle (2pt);
	
	\draw[blue] (0,1.5) -- (0,.9);
	
	\node at (0,-1.4) {$\textbf{.}$};
	
	\node at (0,1.5) {$\textbf{.}$};
	
	\node at (0,.9) {$\textbf{.}$};

	\draw[blue] (0,-.8) arc (90:450: 2.3 and -.85);
	
	\draw[blue] (0,-1.4) arc (90:450: 2.92 and -1.45);

	\node at (0,-2.5) {$\textbf{
			6-6}$};
	
	\node at (-5,-0.74) {$\textbf{
		}$} ;
	
	\end{tikzpicture}
\end{center}
In the meantime, Cummins has determined in \cite{c1} all the genus 1 congruence subgroups of $\PR$ that are commensurable with $\SL$. From his tables, we can conclude that our groups at index 6 and 12 are all congruence subgroups of $\SL$ with the level being the highest cusp width. In particular, the group of index 6 is the commutator group of $\SL$ and has level 6, while the group of level 11 at index 12 is simply $\G_0(11)$. It is worth mentioning that among the genus 1 torsion-free modular subgroups, only 48 of them are congruence subgroups with the highest index being 108, see \cite{c2} for more details.

In a forthcoming work, more analysis and enumeration of these groups will be carried out. In particular, we study the elliptic curves that arise from their modular curves.

\ 

 {\bf Acknowledgment.} We wish to thank   Chris Cummins and Yifan Yang for helpful discussions

\end{document}